\documentclass[twocolumn]{autart}    % Enable this line and disable the 
\usepackage{graphicx}          % Include this line if your 
\usepackage{amsmath}
\usepackage{amsfonts}
\usepackage{units}
\usepackage{csquotes}
\usepackage{color}
\usepackage{mathptmx}

\allowdisplaybreaks

\begin{document}

\begin{frontmatter}

\title{Stop-and-Go Suppression in Two-Class Congested Traffic}            

\thanks[*]{Corresponding author.}

\author[Stuttgart]{Mark Burkhardt}\ead{mburkhardt@eng.ucsd.edu},    % Add the 
\author[UCSD,*]{Huan Yu}\ead{huy015@ucsd.edu}  \text{ }and               % e-mail address 
\author[UCSD]{Miroslav Krstic}\ead{krstic@ucsd.edu}  % (ead) as shown

\address[Stuttgart]{Institute for System Dynamics, University of Stuttgart, Waldburgstr. 17/19, Stuttgart, 70563, Germany}             
\address[UCSD]{Department of Mechanical and Aerospace Engineering, University of California, San Diego, La Jolla, CA, 92093, USA}  % Please supply                                              
% full addresses here.

                             % keyword list or with the 
                                          % help of the Automatica 
                                          % keyword wizard

\begin{abstract}                          % Abstract of not more than 200 words.
This paper develops boundary feedback control laws in order to damp out traffic oscillations in the congested regime of the linearized two-class Aw-Rascle (AR) traffic model. The macroscopic second-order two-class AR traffic model consists of four hyperbolic partial differential equations (PDEs) describing the dynamics of densities and velocities on freeway. The concept of area occupancy is used to express the traffic pressure and equilibrium speed relationship yielding a coupling between the two classes of vehicles. Each vehicle class is characterized by its own vehicle size and driver's behavior. The considered equilibrium profiles of the model represent evenly distributed traffic with constant densities and velocities of both classes along the investigated track section. After linearizing the model equations around those equilibrium profiles, it is observed that in the congested traffic one of the four characteristic speeds is negative, whereas the remaining three are positive. Backstepping control design is employed to stabilize the $4 \times 4$ heterodirectional hyperbolic PDEs. The control input actuates the traffic flow at outlet of the investigated track section and is realized by a ramp metering. A full-state feedback is designed to achieve finite time convergence of the density and velocity perturbations to the equilibrium at zero. This result is then combined with an anti-collocated observer design in order to construct an output feedback control law that damps out stop-and-go waves in finite time by measuring the velocities and densities of both vehicle classes at the inlet of the investigated track section. The performance of the developed controllers is verified by simulation.
\end{abstract}

\end{frontmatter}

\section{Introduction}

Nowadays, more and more people own a car leading to crowded highways and congested traffic during rush hours in many countries of the world. Stop-and-go traffic is common to appear in congested traffic. This phenonemon is characterized by traffic density and velocity perturbations, causing higher fuel consumption and a higher risk of accidents. Traffic management systems like ramp metering or variable speed limits exist in order to regulate traffic. Moreover, they can be used to damp out those traffic oscillations and distribute traffic evenly on the highway. Many recent efforts are focused on boundary control of stop-and-go traffic including~\cite{ARZYU},~\cite{ZhangPrieur_ExpStability_Positivehyp},~\cite{VSLARZ} and~\cite{2Lane_Yu}. However, the results presented in this paper assume heterogeneous vehicle sizes and drivers' behavior. Thus, the overall challenge addressed in this work is the design of a ramp metering traffic management system to reduce traffic oscillations in the congested regime while distinguishing two different vehicle classes. 

In general, traffic models are categorized in micro-, meso- and macroscopic models. Macroscopic models describe the traffic as a continuum and are thus more suitable to investigate traffic jams and disturbances in traffic flow. Typically, their model equations are PDEs. The first representative of this model category is the Lighthill-Whitham-Richards (LWR) model,~\cite{LWR_model1} and~\cite{LWR_model2}, that is given by a single PDE conserving the traffic flow. Although capturing a good amount of realistic traffic phenomena, the LWR model fails to model important phenomena like platoon dispersion~\cite{extendedLWR} or stop-and-go waves and is faulty under light traffic conditions. Thus, a second order extension is provided independently by~\cite{Payne_PWmodel} and~\cite{Whitham_PWmodel}. This second order extension is denoted as the PW model and introduces a second PDE modeling the velocity dynamics of the vehicles in addition to the conversation law. However, the critique in~\cite{DaganzoRequiem} points out that the PW model fails to portray that personalities of vehicles remain unchanged and drivers react more likely to events in front of them. To overcome these issues,~\cite{AwRascleResur} and~\cite{ARZ_Zhang} presented the Aw-Rascle-Zhang (ARZ) model which is also a second-order model describing the interaction between the vehicles by a traffic pressure function. 

Macroscopic multi-class models were proposed after the publication of the ARZ model, including the first order models~\cite{npopulations}, \cite{FOMC_3DFD}, \cite{creeping}, \cite{Logghe_PhD}, \cite{FOMC_porousflow}, \cite{FOMC_Ngoduy}, \cite{FOMC_Fastlane}, \cite{extendedLWR} and \cite{FOMC_Zhang}. The extended LWR model, introduced by~\cite{extendedLWR}, is the first macroscopic multi-class extension. Instead of adjusting their velocity only depending on the density of their own class, vehicles are affected by the densities of all other classes. Thus, the assumed speed-density relationship is formulated with respect to the sum of all densities, the total density. Furthermore, the $n$-populations model~\cite{npopulations} extends this idea of coupling regarding the average length of the vehicle classes by denoting the speed-density relationship in dependence of the mean free space between the vehicles. In~\cite{creeping} and~\cite{FOMC_porousflow}, a phenomenon called creeping is introduced. This behavior occurs in reality and corresponds to the scenario where one vehicle class, for instance trucks, are jammed and stopped due to congestion and a second vehicle class, for instance motorcycles, still moves in the gap between the trucks. 

In addition to first-order macroscopic multi-class models, second-order multi-class models are introduced in~\cite{DCMCLWR}, \cite{SOMC_Gupta}, \cite{SOMC_ExtendedSG}, \cite{MCAR}, \cite{SOMC_TangAnisotropic} and \cite{SOMC_TangDerviedfromCarfollow}. While first-order models assume that the velocities of all vehicles equal to their equilibrium velocities at every time, second-order models provide PDEs describing the velocity dynamics for each class. The denoted second-order models differ in the terms that occur in the velocity dynamics and in the principles which are used to deduce them. For instance,~\cite{SOMC_TangDerviedfromCarfollow} presents an originally microscopic model that is transformed to the macroscopic one, whereas the model in~\cite{MCAR} is based on the macroscopic second-order model of~\cite{AwRascleResur}. The main focus of this paper is the macroscopic multi-class model~\cite{MCAR}. Considering this extension for the case of two different classes yields four coupled nonlinear hyperbolic PDEs which are denoted as the two-class AR traffic model in the following. In order to consider vehicle sizes, the vehicles are assumed to adjust their speed according to a measure called area occupancy,~\cite{AO_Definition_solala} and~\cite{AO_Definition_gut}, which needs to be distinguished from occupancy~\cite{O_Definition}.

Typically, traffic management systems act on the boundary of the investigated track section yielding a boundary control problem. In the literature, different techniques are proposed that achieve convergence of the states of hyperbolic coupled PDEs to a constant equilibrium with boundary control. The main focus of this paper is on the backstepping stabilization technique. While~\cite{Deutscher_2x2_1}, \cite{Deutscher_2x2_2}, \cite{DiMeglio_3x3Backstepping}, \cite{Hu_FiniteTime3x3}, \cite{JMC_H2_Stabilization}, \cite{Vazquez_2x2_CntrlObserv} provide results in case of $2\times 2$ or $3 \times 3$ coupled hyperbolic systems, there also exists literature if an arbitrary amount of linear coupled PDEs is considered~\cite{AuriolControl_nxm}, \cite{Deutscher_nxm}, \cite{DiMeglio_Control}, \cite{HuKrsticGeneralPDEs}, \cite{Su_Nonlocal}. In fact, the presented output feedback control of this work corresponds to the special case of the theoretical result in~\cite{HuKrsticGeneralPDEs} for $m=3$ and $n=1$. The elimination of traffic oscillations by applying the backstepping technique is achieved for the ARZ model in~\cite{VSLARZ} and in~\cite{ARZYU}. Moreover,~\cite{ARZYU} also presents adaptive control results.

Main contribution of this paper: this work presents the first result on boundary control design for traffic congestion consisting of two different vehicle classes. On one hand, this work contributes to traffic modeling in the sense of deducing a macroscopic multi-class traffic model in its characteristic form and investigating the obtained characteristic speeds. On the other hand, a connection between the theoretical control design method backstepping and an up-to-date
extension of the AR traffic model for two classes is created by designing a full-state feedback controller and output feedback controller as ramp metering signal. Moreover, the results are a first step for control problem of more than two vehicle classes or the combination of multiple classes with multiple lanes. 

The rest of this paper is structured as follows: Section $2$ introduces the two-class AR traffic model, the parameters characterizing the two classes and where they occur as well as the assumed boundary conditions. Section $3$ includes the preparation of the linearized model for the control design and the formulation of the control design model. Furthermore, the full-state feedback controller result is presented in section $4$ and the following section $5$ presents the anti-collocated observer design. Section $6$ provides the combination of the results to derive an the output feedback controller and section $7$ verifies the performance of the presented controllers with simulation results. Some future work is discussed in section~$8$.

\section{Problem statement}
%In a first step, the Two Class Aw-Rascle Traffic Model is presented and important model parameters are explained. Afterwards, the model equations are linearized around a constant equilibrium, followed by a discussion of the control objective as well as the control input that is assumed to be available. Finally, it is investigated for which conditions the phenomena, which are about to be eliminated by control, occur.
The two-class AR traffic model is presented and important model parameters are explained. Afterwards, the model equations are linearized around a constant equilibrium, followed by a discussion of the free-flow or congested regimes.
\subsection{Two-class AR traffic model}
The Extended AR model for heterogeneous traffic presented in~\cite{MCAR} is investigated in case of two classes. This two-class AR traffic model is given by
\begin{align}
\partial_t \rho_1  =& -\partial_x(\rho_1 v_1), \label{sec2:eq:TwoClassAR_nonlin1} \\
\partial_t (v_1+p_1(AO)) =& -v_1\partial_x (v_1+p_1(AO)) \notag \\
&+\frac{V_{e,1}(AO)-v_1}{\tau_1} , \label{sec2:eq:TwoClassAR_nonlin2} \\
\partial_t \rho_2  =& -\partial_x(\rho_2 v_2), \label{sec2:eq:TwoClassAR_nonlin3}\\
\partial_t (v_2+p_2(AO)) =& -v_2\partial_x (v_2+p_2(AO)) \notag \\
&+\frac{V_{e,2}(AO)-v_2}{\tau_2}, \label{sec2:eq:TwoClassAR_nonlin4}
\end{align}
where each vehicle class is described by traffic density $\rho_i(x,t)$ and velocity $v_i(x,t)$ with $(x,t) \in [0,L]\times [0,\infty)$. The parameter $L$ is the length of the investigated track section. The traffic density $\rho_i(x,t)$ is defined as vehicles per unit length. The higher the traffic density, the more crowded is the traffic of class $i$ vehicles at a specific spatial point. In addition, the velocity $v_i(x,t)$ describes the velocity of class $i$ vehicles at a specified spatial point along the investigated track section. The traffic density $\rho_1(x,t)$ and velocity $v_1(x,t)$ correspond to the first vehicle class and the traffic density $\rho_2(x,t)$ and velocity $v_2(x,t)$ correspond to the second vehicle class. The non-zero terms on the right hand side represent the adaption of the vehicles to their desired velocities, where $\tau_i$ is the adaptation time.

\begin{figure*}[htbp]
\begin{center}
% trim: links, unten, rechts, oben
\includegraphics[width=16cm,height=7cm, trim = {0cm 1cm 0cm 1cm},clip]{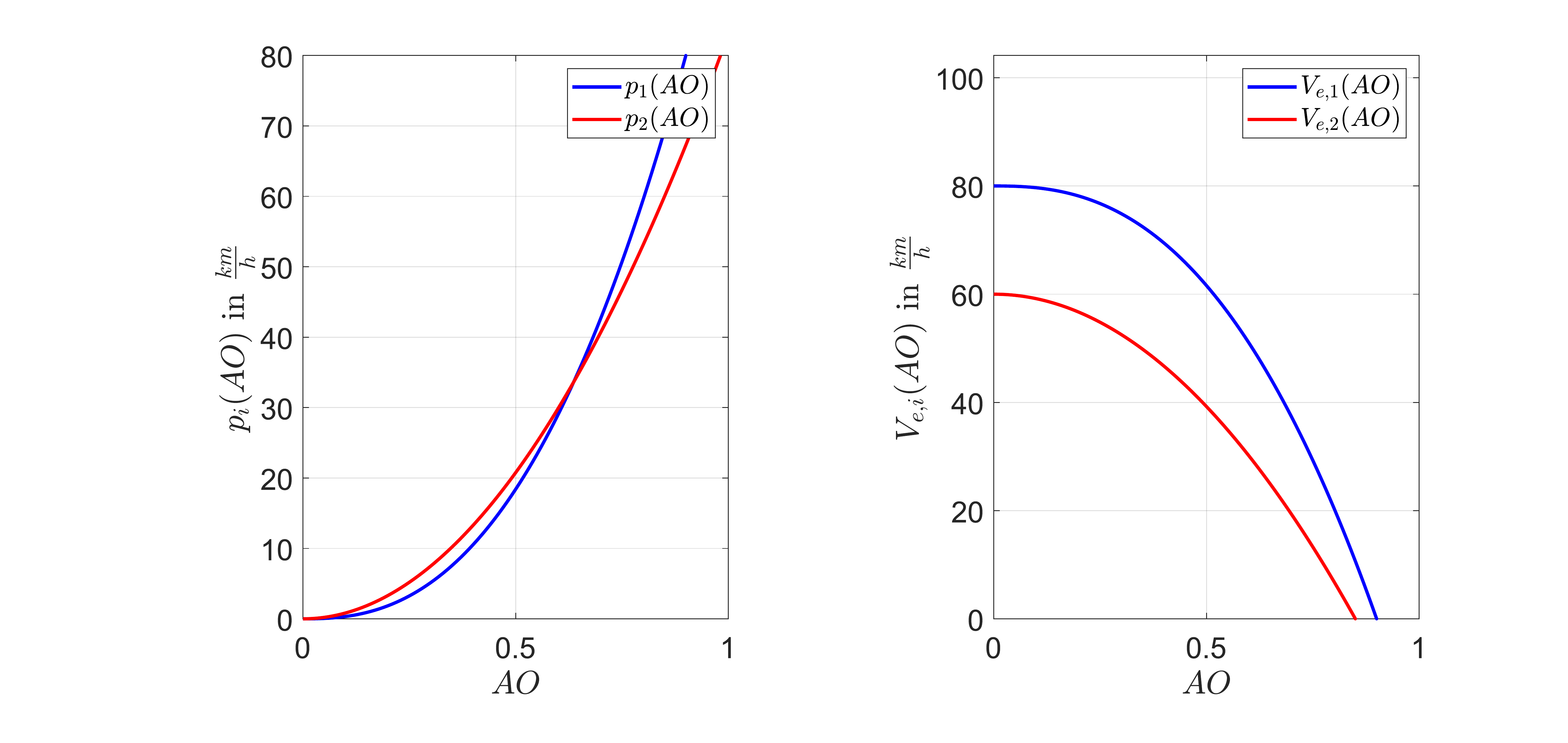}    % The printed column  
\caption{Traffic pressure functions $p_1(AO)$ and $p_2(AO)$ (left) and equilibrium speed-$AO$ relationships $V_{e,1}(AO)$ and $V_{e,2}(AO)$ (right) for the example parameter set $\gamma_1 = 2.5$, $V_1 = \unit[80]{\frac{km}{h}}$, $\overline{AO}=0.9$ for class $1$ and $\gamma_2=2$, $V_2=\unit[60]{\frac{km}{h}}$, $\overline{AO}_2=0.85$ for class $2$.}	% width is 8.4 cm.
\label{sec2:fig:ExamplePlot_Pressure_EquilSpeed}                              % Size the figures 
\end{center}                                 % accordingly.
\end{figure*}

The variable $AO(\rho_1,\rho_2)$ describes the area occupancy and is based on the definition of the occupancy introduced in~\cite{O_Definition}. In~\cite{MCAR}, the expression for the area occupancy is simplified to
\begin{equation}
\label{sec2:eq:AO_def}
AO(\rho_1,\rho_2)=\frac{a_1L\rho_1+a_2L\rho_2}{WL},
\end{equation}
where $a_i$ is the occupied surface per vehicle class $i$ and $W$ the width of the investigated track. Assuming that the traffic densities are $\rho_1(x,t)$ and $\rho_2(x,t)$ along the entire considered highway section, the area occupancy $AO(\rho_1,\rho_2)$ is the percentage of occupied road space by \textit{any} class. It holds that $0\leq AO\leq 1$. The area occupancy depends on both densities since the occupied road surface is influenced by the vehicles of both classes. \\
The traffic pressure function $p_i(AO)$ is formulated as
\begin{equation}
\label{sec2:eq:Pressure_AO}
p_i(AO) = V_{i}\left(\frac{AO(\rho_1,\rho_2)}{\overline{AO}_{i}}\right)^{\gamma_i},
\end{equation}
where $V_i$ corresponds to the free-flow velocity, $\gamma_i>1$ to the traffic pressure exponent and $0<\overline{AO}_i\leq 1$ to the maximum area occupancy. The traffic pressure $p_i(AO)$ is the experienced traffic pressure by class $i$ vehicles and depends on the area occupancy. The higher the area occupancy, the higher the experienced traffic pressure. For instance, if a vehicle suddenly decelerates, then the following vehicle experiences a high traffic pressure forcing another deceleration. Thereby, the free-flow velocity $V_i$ represents the desired velocity of a driver, if no other vehicles of any class are present. The pressure exponent $\gamma_i$ is a parameter that models the experience of the traffic pressure. Higher traffic pressure exponents lead to less experienced pressure. However, the maximum experienced pressure remains the same and is given by the free-flow velocity. The maximum area occupancy $\overline{AO}_i$ describes the percentage of occupied road surface for which the corresponding vehicle class is jammed. To obtain physically meaningful results, $0 < \overline{AO}\leq 1$ holds. For instance, $\overline{AO}_2=0.8$ means that if $80\%$ of the highway are covered by vehicles of \textit{any} class, then the class $2$ vehicles are jammed and therefore their desired velocity is zero. Finally, the equilibrium speed-$AO$ relationship is
\begin{equation}
\label{sec2:eq:Equil_speed_AO}
V_{e,i}(AO) = V_{i}\left(1-\left(\frac{AO(\rho_1,\rho_2)}{\overline{AO}_{i}}\right)^{\gamma_i}\right),
\end{equation}
according to the model of Greenshield \cite{Greenshield}, and represents the desired velocity of the class $i$ vehicles. It depends on the area occupancy since a very crowded road implies a lower desired speed in contrast to a nearly empty road. If the area occupancy is at the maximum $\overline{AO}_i$, then the corresponding equilibrium speed-$AO$ relationship value is $V_{e,i}(\overline{AO}_i) = 0$. In order to show the qualitative behavior of the traffic pressure function~\eqref{sec2:eq:Pressure_AO} and the equilibrium speed-$AO$ relationship~\eqref{sec2:eq:Equil_speed_AO}, both functions are plotted in Figure~\ref{sec2:fig:ExamplePlot_Pressure_EquilSpeed} using an example parameter set. It is illustrated, that a more crowded highway, corresponding to a higher area occupancy, implies a higher experienced traffic pressure and a lower equilibrium speed.
\subsection{Linearized two-class AR traffic model}
The two-class AR traffic model~\eqref{sec2:eq:TwoClassAR_nonlin1} to~\eqref{sec2:eq:TwoClassAR_nonlin4} is linearized around a constant equilibrium state $(\rho_1^*,v_1^*,\rho_2^*,v_2^*)^T$. Inserting this constant state in~\eqref{sec2:eq:TwoClassAR_nonlin1} to~\eqref{sec2:eq:TwoClassAR_nonlin4} yields the conditions
\begin{align}
v_1^*(\rho_1^*,\rho_2^*) &= V_{e,1}(AO(\rho_1^*,\rho_2^*)), \label{sec2:eq:EquilibriumState_cond1} \\
v_2^*(\rho_1^*,\rho_2^*) &= V_{e,2}(AO(\rho_1^*,\rho_2^*)). \label{sec2:eq:EquilibriumState_cond2}
\end{align}
Thus, the equilibrium velocities are determined by the equilibrium densities $\rho_1^*$ and $\rho_2^*$. The perturbations of the distributed variables $\rho_i(x,t)$ and $v_i(x,t)$ are defined as
\begin{align}
\tilde{\rho}_i(x,t) = \rho_i(x,t)-\rho_i^*, \label{sec2:eq:PerturbationsDefinition1} \\
\tilde{v}_i(x,t) = v_i(x,t)-v_i^*, \label{sec2:eq:PerturbationsDefinition2}
\end{align}
for each class $i$ and the linearized model equations are given by
\begin{equation}
\label{sec2:eq:LinearizedModelEq}
J_t\left[\begin{array}{c}
\tilde{\rho}_{1t} \\
\tilde{v}_{1t} \\
\tilde{\rho}_{2t} \\
\tilde{v}_{2t} 
\end{array}\right]+J_x\left[\begin{array}{c}
\tilde{\rho}_{1x} \\
\tilde{v}_{1x} \\
\tilde{\rho}_{2x} \\
\tilde{v}_{2x} \\
\end{array}\right] + J\left[\begin{array}{c}
\tilde{\rho}_1 \\
\tilde{v}_1 \\
\tilde{\rho}_2 \\
\tilde{v}_2
\end{array}\right] = \left[\begin{array}{c}
0 \\
0 \\
0 \\
0 \\
\end{array}\right],
\end{equation}
where the introduced Jacobian matrices are  
\begin{align}
J_t &= \left[\begin{array}{cccc}
1 & 0 & 0 & 0 \\
\beta_{11} & 1 & \beta_{12} & 0 \\
0 & 0 & 1 & 0 \\
\beta_{21} & 0 & \beta_{22} & 1 \\
\end{array}\right], \\
J_x &= \left[\begin{array}{cccc}
v_1^* & \rho_1^* & 0 & 0 \\
v_1^*\beta_{11} & v_1^* & v_1^* \beta_{12} & 0 \\
0 & 0 & v_2^* & \rho_2^* \\
v_2^*\beta_{21} & 0 & v_2^*\beta_{22} & v_2^* \\
\end{array}\right], \\
J &= \left[\begin{array}{cccc}
0 & 0 & 0 & 0 \\
\frac{1}{\tau_1}\beta_{11} & \frac{1}{\tau_1} & \frac{1}{\tau_1} \beta_{12} & 0 \\
0 & 0 & 0 & 0 \\
\frac{1}{\tau_2}\beta_{21} & 0 & \frac{1}{\tau_2}\beta_{22}   & \frac{1}{\tau_2} \\
\end{array}\right],
\end{align}
and the abbreviations 
\begin{equation}
\beta_{ij}(\rho_1^*,\rho_2^*) = \left.\frac{\partial p_i(AO(\rho_1,\rho_2))}{\partial \rho_j}\right|_{\rho_1= \rho_1^*,\rho_2 = \rho_2^*}
\end{equation}
are introduced with $i,j=1,2$. The abbreviations $\beta_{ij}(\rho_1^*,\rho_2^*)$ represent the derivative of the class $i$ traffic pressure function with respect to class $j$ traffic density. The boundary conditions are assumed to be
\begin{align}
\rho_1(0,t) &= \rho_1^*,\label{sec2:eq:BC0_nonlin_Dens1} \\
\rho_2(0,t) &= \rho_2^*, \label{sec2:eq:BC0_nonlin_Dens2} \\
q_1(0,t)+ q_2(0,t) &= q_1^*+q_2^*, \label{sec2:eq:BC0_nonlin_Flow} \\
q_1(L,t)+q_2(L,t) &= q_1^*+q_2^*, \label{sec2:eq:BCL_nonlin_Flow_woContrl}
\end{align}
where~\eqref{sec2:eq:BC0_nonlin_Flow} and~\eqref{sec2:eq:BCL_nonlin_Flow_woContrl} assume that the same total traffic flow enters and leaves the track section which is given by the sum of the class $1$ and class $2$ equilibrium flows $q_1^*$ and $q_2^*$. The traffic flow of class $i$ is defined as
\begin{equation}
q_i(x,t) = \rho_i(x,t)v_i(x,t).
\end{equation}
Boundary conditions~\eqref{sec2:eq:BC0_nonlin_Dens1} and~\eqref{sec2:eq:BC0_nonlin_Dens2} indicate that the traffic densities of the incoming traffic flow are equivalent to the equilibrium densities. Thus, not only the entering traffic flow is constant, in fact the densities of both classes in this traffic flow are assumed to be constant. The linearization of the introduced boundary conditions~\eqref{sec2:eq:BC0_nonlin_Dens1} to~\eqref{sec2:eq:BCL_nonlin_Flow_woContrl} is
\begin{align}
0 &= \tilde{\rho}_1(0,t), \label{sec2:eq:BC0_lin_Dens1} \\
0 &= \tilde{\rho}_2(0,t), \label{sec2:eq:BC0_lin_Dens2}\\
0 &= v_1^*\tilde{\rho}_1(0,t)+\rho_1^*\tilde{v}_1(0,t)+v_2^*\tilde{\rho}_2(0,t)+\rho_2^*\tilde{v}_2(0,t),  \label{sec2:eq:BC0_lin_Flow}\\
0 &= v_1^*\tilde{\rho}_1(L,t)+\rho_1^*\tilde{v}_1(L,t)+v_2^*\tilde{\rho}_2(L,t)+\rho_2^*\tilde{v}_2(L,t). \label{sec2:eq:BCL_lin_Flow_woContrl} 
\end{align}
\subsection{Free/congested regime analysis}
In general, two different regimes of traffic are distinguished: the free-flow regime and the congested regime. The free-flow regime is characterized by the fact that the total information of the system travels downstream along with the vehicles. In that case, the model equations correspond to four homodirectional hyperbolic PDEs. On the other hand, a partial upstream propagation of information characterizes the traffic flow in the congested regime. The corresponding heterodirectional behavior causes the development of stop-and-go traffic which implies increased fuel consumption and risk of accidents. Therefore, it is reasonable to investigate which choices of equilibrium densities and parameters lead to heterodirectional information propagation. Therefore, the characteristic speeds are computed and their signs are considered in the following. First, the linearized model equations~\eqref{sec2:eq:LinearizedModelEq} need to be decoupled in time leading to 
\begin{equation}
\left[\begin{array}{c}
\tilde{\rho}_{1t} \\
\tilde{v}_{1t} \\
\tilde{\rho}_{2t} \\
\tilde{v}_{2t}
\end{array}\right] + \tilde{J}_x\left[\begin{array}{c}
\tilde{\rho}_{1x} \\
\tilde{v}_{1x} \\
\tilde{\rho}_{2x} \\
\tilde{v}_{2x} \\
\end{array}\right]=\tilde{J}\left[\begin{array}{c}
\tilde{\rho}_1 \\
\tilde{v}_1 \\
\tilde{\rho}_2 \\
\tilde{v}_2 \\
\end{array}\right]
\end{equation}
with the new Jacobian matrices
\begin{align}
&\tilde{J}_x = \\
&\left[\begin{array}{cccc}
v_1^* & \rho_1^* & 0 & 0 \\
0 & v_1^*-\beta_{11}\rho_1^* & \beta_{12}(v_1^*-v_2^*) & -\beta_{12}\rho_2^* \\
0 & 0 & v_2^* & \rho_2^* \\
\beta_{21}(v_2^*-v_1^*) & -\beta_{21}\rho_1^* & 0 & v_2^*-\beta_{22}\rho_2^* \\
\end{array}\right]\notag
\end{align}
and
\begin{equation}
\tilde{J} = \left[\begin{array}{cccc}
0 & 0 & 0 & 0 \\
-\frac{1}{\tau_1}\beta_{11} & -\frac{1}{\tau_1} & -\frac{1}{\tau_1}\beta_{12} & 0 \\
0 & 0 & 0 & 0 \\
-\frac{1}{\tau_2}\beta_{21} & 0 & -\frac{1}{\tau_2}\beta_{22} & -\frac{1}{\tau_2} \\
\end{array}\right].
\end{equation}
Then, the characteristic speeds are given by the eigenvalues of $\tilde{J}_x$ which are
\begin{align}
\lambda_1 =& v_1^*(\rho_1^*,\rho_2^*), \label{sec2:eq:CharSpeed1} \\
\lambda_2 =& v_2^*(\rho_1^*,\rho_2^*), \label{sec2:eq:CharSpeed2} \\
\lambda_3 =& \frac{v_1^*(\rho_1^*,\rho_2^*)+v_2^*(\rho_1^*,\rho_2^*)-\alpha_1(\rho_1^*,\rho_2^*)-\alpha_2(\rho_1^*,\rho_2^*)}{2} \notag\\
&+\frac{\Delta(\rho_1^*,\rho_2^*)}{2}, \label{sec2:eq:CharSpeed3}\\
\lambda_4 =& \frac{v_1^*(\rho_1^*,\rho_2^*)+v_2^*(\rho_1^*,\rho_2^*)-\alpha_1(\rho_1^*,\rho_2^*)-\alpha_2(\rho_1^*,\rho_2^*)}{2} \notag \\
&-\frac{\Delta(\rho_1^*,\rho_2^*)}{2},  \label{sec2:eq:CharSpeed4}
\end{align}
where 
\begin{equation}
\Delta(\rho_1^*,\rho_2^*) = \sqrt{\left(\alpha_2\rho_2^*-\alpha_1\rho_1^*+v_1^*-v_2^*\right)^2+4\alpha_1\alpha_2\rho_1^*\rho_2^*}
\end{equation}
and 
\begin{equation}
\alpha_i(\rho_1^*,\rho_2^*) = \beta_{ii}(\rho_1^*,\rho_2^*) = \left.\frac{\partial p_i(AO(\rho_1,\rho_2))}{\partial \rho_i}\right|_{\rho_1= \rho_1^*,\rho_2 = \rho_2^*}.
\end{equation}
For model validity, the equilibrium velocities of both vehicle classes are chosen to be positive, i.e. $v_1^*>0$ and $v_2^*>0$ and all vehicles travel downstream. Thus, the first two characteristic speeds~\eqref{sec2:eq:CharSpeed1} and~\eqref{sec2:eq:CharSpeed2} are positive. In addition, it is shown in~\cite{AnalyseARZ_CharSpeed} that
\begin{equation}
\label{sec2:eq:CharSpeed_Order}
\lambda_4 \leq \min\{\lambda_1,\lambda_2\} \leq \lambda_3 \leq \max\{\lambda_1,\lambda_2\} 
\end{equation}
holds. Because $\lambda_1>0$ and $\lambda_2>0$,~\eqref{sec2:eq:CharSpeed_Order} implies that $\lambda_3$ is positive as well. Hence, the only characteristic speed that may have a negative sign is $\lambda_4$. Therefore, traffic is defined to be in the free regime if the equilibrium densities and parameters satisfy
\begin{equation}
\label{sec2:FreeFlowReg_Def}
\lambda_1,\lambda_2,\lambda_3,\lambda_4 > 0
\end{equation}
and in the congested regime if they meet
\begin{equation}
\label{sec2:CongestedReg_Def}
\lambda_1,\lambda_2,\lambda_3 > 0, \quad \lambda_4 <0.
\end{equation}
Since a controller dealing with congested traffic is designed later on, it is assumed that the equilibrium densities and parameters are chosen such that~\eqref{sec2:CongestedReg_Def} holds throughout the rest of this paper. In fact, this result allows to formulate a traffic froude number as it is presented in~\cite{Froude_Belletti} for the ARZ model. Besides, in~\cite{AnalyseARZ_CharSpeed}, a physical interpretation of the characteristic speeds~\eqref{sec2:eq:CharSpeed1} to~\eqref{sec2:eq:CharSpeed4} is given. While $\lambda_1$ and $\lambda_2$ correspond to the flow of class $1$ vehicles and class $2$ vehicles, $\lambda_3$ is related to the fact that the faster vehicle class overtakes the slower one. For that reason, it is reasonable to obtain $\lambda_1 =\lambda_2=\lambda_3$, if $v_1^*=v_2^*$ is assumed.
\begin{figure*}[htbp]
\begin{center}
% trim: links, unten, rechts, oben
\includegraphics[width=16cm,height=7cm, trim = {0cm 0cm 0cm 0cm},clip]{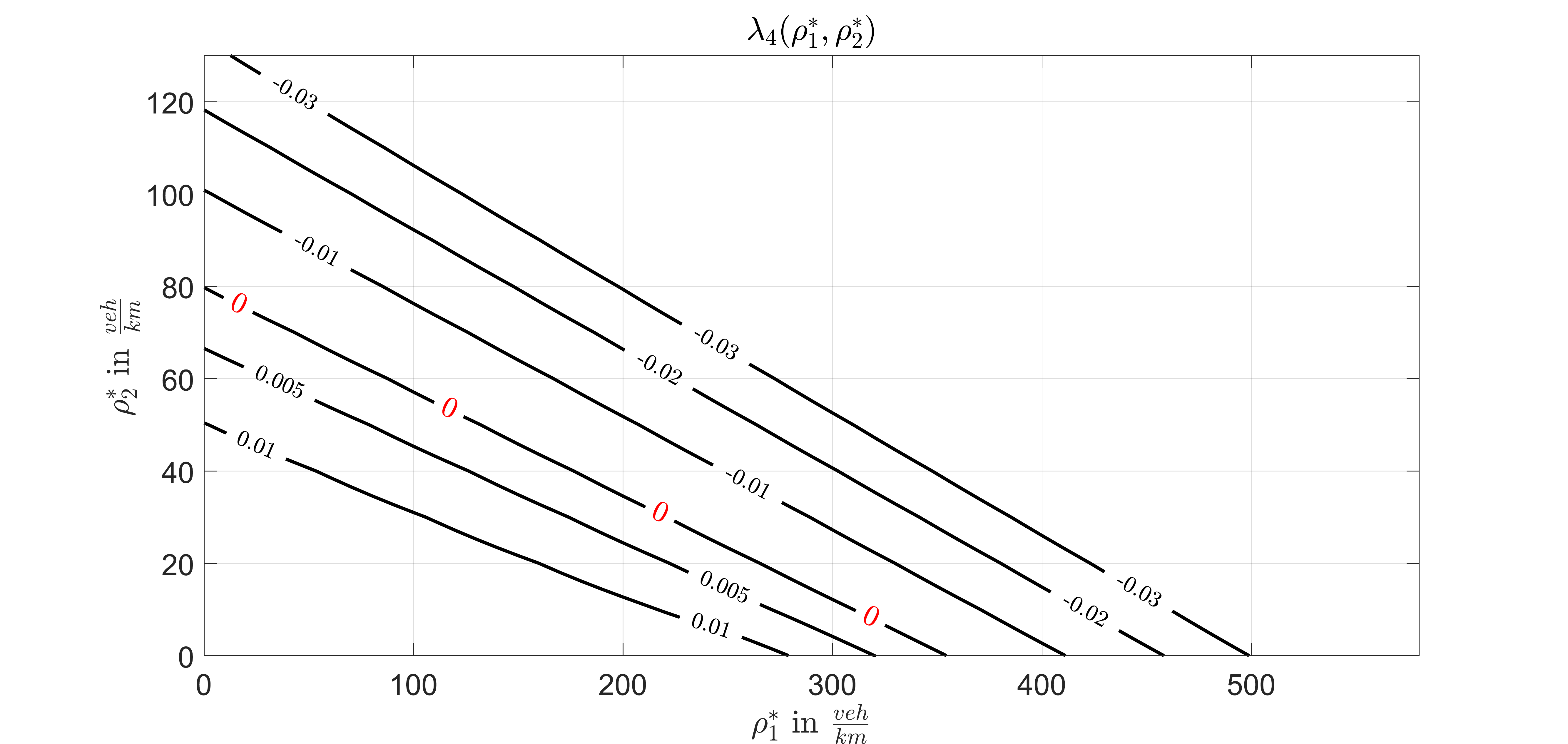}    % The printed column  
\caption{Contour plot of $\lambda_4$ for the parameter set $\gamma_1 = 2.5$, $V_1 = \unit[80]{\frac{km}{h}}$, $\overline{AO}=0.9$ for class $1$ and $\gamma_2=2$, $V_2=\unit[60]{\frac{km}{h}}$, $\overline{AO}_2=0.85$ for class $2$. The contour line $\lambda_4=0$ describes the boundary between the free-flow and congested regime.}	% width is 8.4 cm.
\label{sec2:fig:Contourplot_CongestionBoundary}                              % Size the figures 
\end{center}                                 % accordingly.
\end{figure*}

According to the definitions in~\eqref{sec2:FreeFlowReg_Def} and~\eqref{sec2:CongestedReg_Def}, the boundary between the two regimes is defined as $\lambda_4 = 0$. In case of the two-class AR traffic model, this boundary is a line which can be drawn in the $\rho_1^*$-$\rho_2^*$-plane. Compared to the single class consideration, the boundary is equivalent to a single density which is the critical density. The numerically computed boundary between the two regimes is plotted as a contour plot for an example parameter set in Figure~\ref{sec2:fig:Contourplot_CongestionBoundary}. The figure also indicates that small values for both equilibrium densities $\rho_1^*$ and $\rho_2^*$ correspond to a positive value of $\lambda_4$ and therefore homodirectional behavior. On the other hand, large values of densities lead to a negative value of $\lambda_4$ indicating heterodirectional behavior. Hence, smaller equilibrium densities correspond to free-flow regime and large equilibrium densities correspond to congested regime. \\
%Another interesting difference between the single class and two class consideration, are the densities that maximize the traffic flow. In the single class consideration, the fundamental diagram describes that the maximum traffic flow is achieved at the critical density. However, in case of the Two-Class Aw-Rascle Traffic Model, the densities that maximize the traffic flow are given by the area occupancy value that maximizes the equilibrium flow-$AO$ relationship, defined as
%\begin{align}
%Q_e(AO) &= Q_{e,1}(AO)+Q_{e,2}(AO) \notag\\
%	    &= V_{e,1}(AO)AO+V_{e,2}(AO)AO.
%\end{align}
%Computing this area occupancy value $\check{AO}(\rho_1^*,\rho_2^*)$ such that $Q_e'(\check{AO}(\rho_1^*,\rho_2^*))=0$ yields again a line in the $\rho_1^*$-$\rho_2^*$-plane.
\section{Boundary control design model}
In a next step, the control objective and input is explained and afterwards the introduced linearized two-class AR traffic model,~\eqref{sec2:eq:LinearizedModelEq},~\eqref{sec2:eq:BC0_lin_Dens1} to~\eqref{sec2:eq:BCL_lin_Flow_woContrl}, is prepared for the control design. The preparation is done using two transformations, the transformation to Riemann coordinates and a second transformation to further simplify the equations expressed in Riemann coordinates resulting in the control design model.
\subsection{Control input and objective}
The overall goal is to damp out stop-and-go traffic in the congested regime and achieve convergence to the equilibrium states in a finite time. The term of stop-and-go traffic refers to oscillations of the density and velocity perturbations around their constant equilibrium values along the highway. Moreover, a ramp metering is considered to be installed at the outlet of the investigated track section regulating the traffic outflow. In this work, the ramp metering is used to damp out the introduced stop-and-go waves. Thus, the boundary condition~\eqref{sec2:eq:BCL_lin_Flow_woContrl} becomes
\begin{equation}
U(t) = v_1^*\tilde{\rho}_1(L,t)+\rho_1^*\tilde{v}_1(L,t)+v_2^*\tilde{\rho}_2(L,t)+\rho_2^*\tilde{v}_2(L,t). \label{sec2:eq:BCL_lin_Flow_wContrl} 
\end{equation}
Compared to the application of multi-phase flow in oil pipelines,~\cite{DiMeglio_Slugging}, a ramp metering works as a valve at the end of the pipe to control the outgoing flow.
\subsection{Transformation to Riemann coordinates}
\begin{figure}[tbp]
\begin{center}
% trim: links, unten, rechts, oben
\includegraphics[width=8.4cm,height=6cm, trim = {5.5cm 4cm 5.5cm 3.5cm},clip]{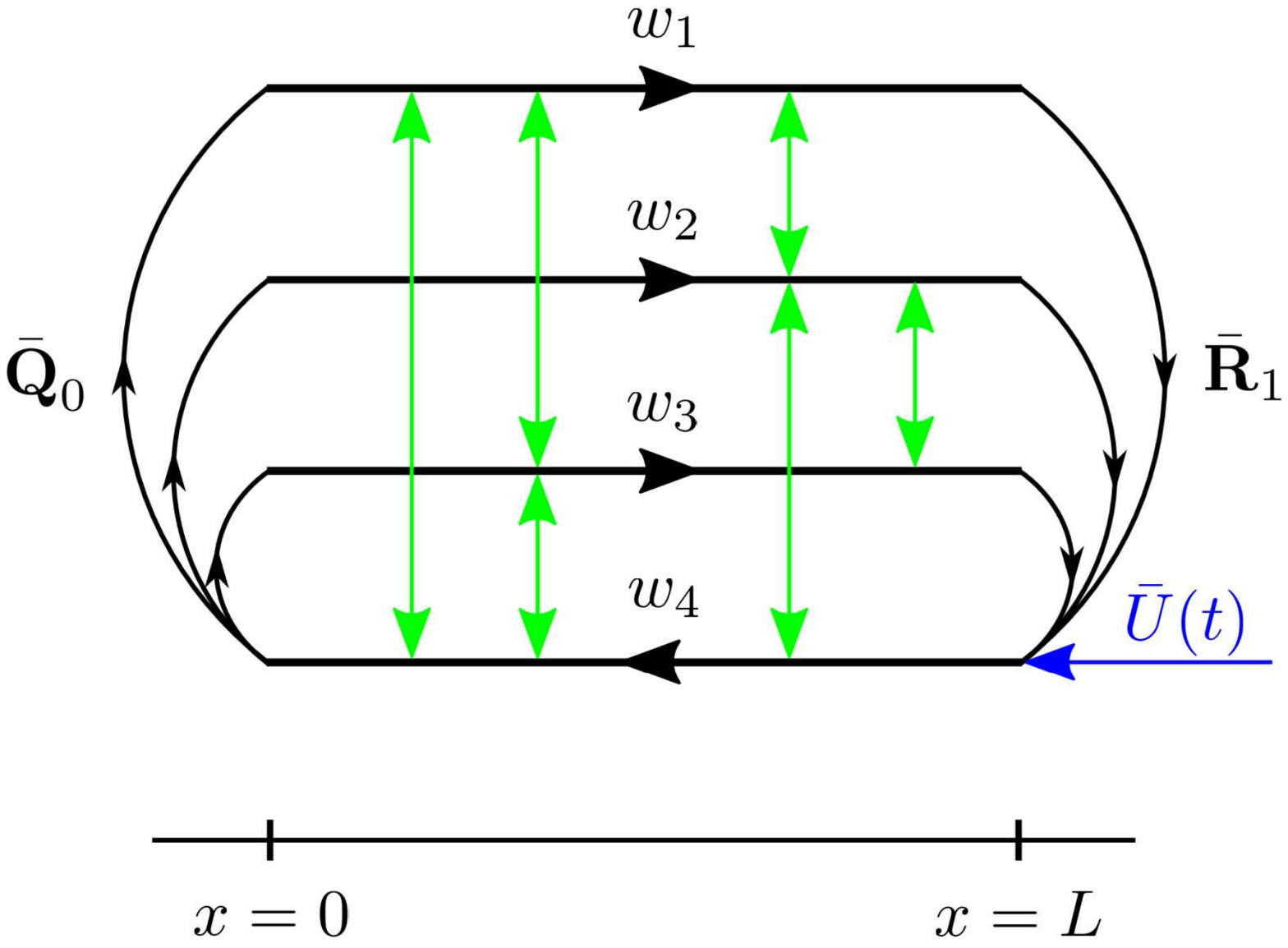}    % The printed column  
\caption{Schematic diagram of the control design model. The green arrow indicates the location where the control input acts on the system. The blue arrows represent the couplings between all four states.}	% width is 8.4 cm.
\label{sec3:fig:Qualitative_ControlDesignModel}                              % Size the figures 
\end{center}                                 % accordingly.
\end{figure}

The system is transformed to Riemann coordinates to accomplish a decoupling of the spatial derivatives. The Riemann variables $(\bar{w}_1,\bar{w}_2,\bar{w}_3,\bar{w}_4)^T$ are defined as new coordinates. The linear state transformation is given by
\begin{equation}
\left[\begin{array}{c}
\bar{w}_1 \\
\bar{w}_2 \\
\bar{w}_3 \\
\bar{w}_4 \\
\end{array}\right] = V^{-1}\left[\begin{array}{c}
\tilde{\rho}_1 \\
\tilde{v}_1 \\
\tilde{\rho}_2 \\
\tilde{v}_2 \\
\end{array}\right]
\end{equation}
where the constant invertible transformation matrix $V$ satisfies
\begin{equation}
\left[\begin{array}{cccc}
\lambda_1 & 0 & 0 & 0 \\
0 & \lambda_2 & 0 & 0 \\
0 & 0 & \lambda_3 & 0 \\
0 & 0 & 0 & \lambda_4 \\
\end{array}\right] = V^{-1}\tilde{J}_xV
\end{equation}
and therefore diagonalizes the Jacobian $\tilde{J}_x$. The entries of $V$ are denoted as 
\begin{equation}
V = \left\{v_{ij}\right\}_{1\leq i\leq 4,1\leq j \leq 4}.
\end{equation}
The matrix $V$ is straightforward to obtain but omitted in this paper due to its complexity and length. Inserting the transformation in~\eqref{sec2:eq:LinearizedModelEq} yields the model equations in Riemann coordinates
\begin{equation}
\label{sec3:eq:PDEs_Riemann}
\left[\begin{array}{c}
\bar{w}_{1t} \\
\bar{w}_{2t} \\
\bar{w}_{3t} \\
\bar{w}_{4t} \\
\end{array}\right]+\left[\begin{array}{cccc}
\lambda_1 & 0 & 0 & 0 \\
0 & \lambda_2 & 0 & 0 \\
0 & 0 & \lambda_3 & 0 \\
0 & 0 & 0 & \lambda_4 \\
\end{array}\right]\left[\begin{array}{c}
\bar{w}_{1x} \\
\bar{w}_{2x} \\
\bar{w}_{3x} \\
\bar{w}_{4x} \\
\end{array}\right] = \hat{J}\left[\begin{array}{c}
\bar{w}_1 \\
\bar{w}_2 \\
\bar{w}_3 \\
\bar{w}_4 \\
\end{array}\right],
\end{equation}
where 
\begin{equation}
\hat{J} = V^{-1}\tilde{J}V
\end{equation}
and the entries of the Jacobian $\hat{J}$ are denoted by
\begin{equation}
\hat{J} = \{\hat{J}_{ij}\}_{1 \leq i \leq 4, 1\leq j \leq 4}.
\end{equation}
Since the coefficient matrix of the spatial derivatives is now diagonal, a decoupling in spatial derivatives is achieved. The characteristic speeds~\eqref{sec2:eq:CharSpeed1} to~\eqref{sec2:eq:CharSpeed4} form the diagonal because they are the eigenvalues of $\tilde{J}_x$. In addition, the same transformation is applied to the boundary conditions~\eqref{sec2:eq:BC0_lin_Dens1},~\eqref{sec2:eq:BC0_lin_Dens2},~\eqref{sec2:eq:BC0_lin_Flow} and~\eqref{sec2:eq:BCL_lin_Flow_wContrl} yielding
\begin{align}
\left[\begin{array}{c}
\bar{w}_1(0,t) \\
\bar{w}_2(0,t) \\
\bar{w}_3(0,t) \\
\end{array}\right] &= \hat{Q}_0\bar{w}_4(0,t), \label{sec3:eq:BC0_Riemann} \\
\bar{w}_4(L,t) &= \hat{R}_1\left[\begin{array}{c}
\bar{w}_1(L,t) \\
\bar{w}_2(L,t) \\
\bar{w}_3(L,t) \\
\end{array}\right]+\hat{U}(t). \label{sec3:eq:BCL_Riemann}
\end{align}
The matrices are given by
\begin{align}
\hat{Q}_0 &= -\left[\begin{array}{ccc}
v_{11} & v_{12} & v_{13} \\
v_{31} & v_{32} & v_{33} \\
\kappa_1 & \kappa_2 & \kappa_3 \\
\end{array}\right]^{-1}\left[\begin{array}{c}
v_{14} \\
v_{34} \\
\kappa_4
\end{array}\right], \label{sec3:eq:hat_Q0} \\
\hat{R}_1 &= -\frac{1}{\kappa_4} \left[\begin{array}{ccc}
\kappa_1 & \kappa_2 & \kappa_3
\end{array}\right] \label{sec3:eq:hat_R1} 
\end{align}
and are obtained by formulating the linearized boundary conditions in matrix form, inserting the transformation law to Riemann coordinates and decoupling afterwards. Besides, in~\eqref{sec3:eq:hat_Q0} and~\eqref{sec3:eq:hat_R1}, the abbreviations
\begin{equation}
\label{sec3:eq:kappa_abbrev}
\kappa_{i} = v_1^*v_{1i}+\rho_1^*v_{2i}+v_2^*v_{3i}+\rho_2^*v_{4i}, \quad i = 1,2,3,4, 
\end{equation} 
are inserted. Furthermore, the input transformation, used in~\eqref{sec3:eq:BCL_Riemann}, is
%\begin{equation}
%\hat{U}(t) = \frac{U(t)}{v_1^*v_{14}+\rho_1^*v_{24}+v_2^*v_{34}+\rho_2^*v_{44}}.
%\end{equation}
\begin{equation}
\hat{U}(t) = \frac{1}{\kappa_4}U(t).
\end{equation}
Overall, the model equations in Riemann coordinates are given by~\eqref{sec3:eq:PDEs_Riemann},~\eqref{sec3:eq:BC0_Riemann} and~\eqref{sec3:eq:BCL_Riemann}. Because the inverse linear state transformation exists, the system expressed in Riemann coordinates shares the same stability properties with the original system. 
%\begin{figure*}[htbp]
%\begin{center}
% trim: links, unten, rechts, oben
%\includegraphics[width=16cm,height=9cm, trim = {0cm 4cm 0cm 2cm},clip]{Figures/%QualitativeBehavior_SortedRiemann.pdf}    % The printed column  
%\caption{Schematic diagram of the control design model. The green arrow indicates the location where the control %input acts on the system. The blue arrows represent the couplings between all four states.}	% width is 8.4 cm.
%\label{sec3:fig:Qualitative_ControlDesignModel}                              % Size the figures 
%\end{center}                                 % accordingly.
%\end{figure*} 
In a next step, a second transformation is performed in order to complete the preparation for the control design and obtain the control design model. The resulting coordinates are $(w_1,w_2,w_3,w_4)^T$. The second transformation achieves zero elements on the diagonal of $\hat{J}$ in~\eqref{sec3:eq:PDEs_Riemann} and sorts the positive characteristic speeds~\eqref{sec2:eq:CharSpeed1} to~\eqref{sec2:eq:CharSpeed3} in ascending order on the diagonal of the coefficient matrix of the spatial derivatives. However, an ascending order is only defined uniquely, if it is known whether $\lambda_1>\lambda_2$ or $\lambda_1<\lambda_2$ holds, according to~\eqref{sec2:eq:CharSpeed_Order}. In the following, it is assumed that class $1$ vehicles represent small and fast average vehicles whereas class $2$ describes big trucks which are large and slow. Thus, for the equilibrium velocities $v_1^*>v_2^*$ holds and therefore the ascending order of positive characteristic speeds is $\lambda_2<\lambda_3<\lambda_1$. Hence, the transformation law
\begin{align}
w_1 = e^{-\frac{\hat{J}_{22}}{v_2^*}x}\bar{w}_2, \label{sec3:eq:ZeroDiag_Transform1} \\
w_2 = e^{-\frac{\hat{J}_{33}}{\lambda_3}x}\bar{w}_3, \label{sec3:eq:ZeroDiag_Transform2} \\
w_3 = e^{-\frac{\hat{J}_{11}}{v_1^*}x}\bar{w}_1, \label{sec3:eq:ZeroDiag_Transform3} \\
w_4 = e^{-\frac{\hat{J}_{44}}{\lambda_4}x}\bar{w}_4 \label{sec3:eq:ZeroDiag_Transform4} 
\end{align}
is applied to~\eqref{sec3:eq:PDEs_Riemann} yielding the transformed PDEs
\begin{align}
\left[\begin{array}{c}
w_{1t} \\
w_{2t} \\
w_{3t} \\
\end{array}\right]+\Lambda^+
\left[\begin{array}{c}
w_{1x} \\
w_{2x} \\
w_{3x} \\
\end{array}\right]&= \Sigma^{++}(x)
\left[\begin{array}{c}
w_1 \\
w_2 \\
w_3 \\
\end{array}\right] + \Sigma^{+-}(x)w_4, \label{sec3:eq:ControlDesignModel_w123PDE} \\
w_{4t}- \Lambda^- w_{4x} &= \Sigma^{-+}(x) \left[\begin{array}{c}
w_1 \\
w_2 \\
w_3 \\
\end{array}\right]. \label{sec3:eq:ControlDesignModel_w4PDE}
\end{align}
with
\begin{align}
\Lambda^+ &= \left[\begin{array}{ccc}
v_2^* & 0 & 0 \\
0 & \lambda_3 & 0 \\
0 & 0 & v_1^* \\
\end{array}\right], \\
\Lambda^- &= -\lambda_4, \\
\Sigma^{++}(x) &= \left[\begin{array}{ccc}
0 & \bar{J}_{12}(x) & \bar{J}_{13}(x) \\
\bar{J}_{21}(x) & 0 & \bar{J}_{23}(x) \\
\bar{J}_{31}(x)& \bar{J}_{32}(x) & 0 \\ 
\end{array}\right], \\
\Sigma^{+-}(x) &= \left[\begin{array}{ccc}
\bar{J}_{14}(x) & \bar{J}_{24}(x) & \bar{J}_{34}(x)
\end{array}\right]^T,\\
\Sigma^{-+}(x) &= \left[\begin{array}{ccc}
\bar{J}_{41}(x) & \bar{J}_{42}(x) & \bar{J}_{43}(x)
\end{array}\right].
\end{align}
The abbreviations for the coefficients of the source term, $\bar{J}_{ij}(x), \text{ } i,j=1,2,3,4,$ are:
\begin{align*}
\bar{J}_{12}(x) &= \hat{J}_{23}e^{\left(\frac{\hat{J}_{33}}{\lambda_3}-\frac{\hat{J}_{22}}{v_2^*}\right)x},\bar{J}_{13}(x) = \hat{J}_{21}e^{\left(\frac{\hat{J}_{11}}{v_1^*}-\frac{\hat{J}_{22}}{v_2^*}\right)x}, \\
\bar{J}_{14}(x) &= \hat{J}_{24}e^{\left(\frac{\hat{J}_{44}}{\lambda_4}-\frac{\hat{J}_{22}}{v_2^*}\right)x},\bar{J}_{21}(x) = \hat{J}_{32}e^{\left(\frac{\hat{J}_{22}}{v_2^*}-\frac{\hat{J}_{33}}{\lambda_3}\right)x}, \\
\bar{J}_{23}(x) &= \hat{J}_{31}e^{\left(\frac{\hat{J}_{11}}{v_1^*}-\frac{\hat{J}_{33}}{\lambda_3}\right)x}, \bar{J}_{24}(x) = \hat{J}_{34}e^{\left(\frac{\hat{J}_{44}}{\lambda_4}-\frac{\hat{J}_{33}}{\lambda_3}\right)x}, \\
\bar{J}_{31}(x) &= \hat{J}_{12}e^{\left(\frac{\hat{J}_{22}}{v_2^*}-\frac{\hat{J}_{11}}{v_1^*}\right)x},\bar{J}_{32}(x) = \hat{J}_{13}e^{\left(\frac{\hat{J}_{33}}{\lambda_3}-\frac{\hat{J}_{11}}{v_1^*}\right)x}, \\
\text{ }\bar{J}_{34}(x) &= \hat{J}_{14}e^{\left(\frac{\hat{J}_{44}}{\lambda_4}-\frac{\hat{J}_{11}}{v_1^*}\right)x}, \bar{J}_{41}(x) = \hat{J}_{42}e^{\left(\frac{\hat{J}_{22}}{v_2^*}-\frac{\hat{J}_{44}}{\lambda_4}\right)x}, \\
\text{ }\bar{J}_{42}(x) &= \hat{J}_{43}e^{\left(\frac{\hat{J}_{33}}{\lambda_3}-\frac{\hat{J}_{44}}{\lambda_4}\right)x},\bar{J}_{43}(x) = \hat{J}_{41}e^{\left(\frac{\hat{J}_{11}}{v_1^*}-\frac{\hat{J}_{44}}{\lambda_4}\right)x}. \\
\end{align*}
The diagonal elements of $\Lambda^+$ are sorted in an ascending order and the relations $\lambda_1 = v_1^*$,~\eqref{sec2:eq:CharSpeed1}, and $\lambda_2=v_2^*$,~\eqref{sec2:eq:CharSpeed2}, are inserted. In addition, $\Sigma^{++}(x)$, $\Sigma^{+-}(x)$ and $\Sigma^{-+}(x)$ are depending on the spatial coordinate. Their entries $\bar{J}_{ij}(x)$ are bounded and either positive or negative on the whole domain, depending on the sign of the corresponding $\hat{J}_{ij}$. Applying the transformation~\eqref{sec3:eq:ZeroDiag_Transform1} to~\eqref{sec3:eq:ZeroDiag_Transform4} on the boundary conditions~\eqref{sec3:eq:BC0_Riemann} and~\eqref{sec3:eq:BCL_Riemann} yields
\begin{align}
\left[\begin{array}{c}
w_1(0,t) \\
w_2(0,t) \\
w_3(0,t) \\
\end{array}\right] &= \bar{Q}_0w_4(0,t), \label{sec3:eq:ControlDesignModel_bc0} \\
w_4(L,t) &= \bar{R}_1\left[\begin{array}{c}
w_1(L,t) \\
w_2(L,t) \\
w_3(L,t) \\
\end{array}\right]+\bar{U}(t) \label{sec3:eq:ControlDesignModel_bcL}
\end{align}
with
\begin{align}
\bar{Q}_0 &= \left[\begin{array}{ccc}
0 & 0 & 1 \\
1 & 0 & 0 \\
0 & 1 & 0 
\end{array}\right]^{-1}\hat{Q}_0, \notag \\
\bar{R}_1 &= \hat{R}_1\left[\begin{array}{ccc}
0 & 0 & e^{\left( \frac{\hat{J}_{11}}{v_1^*}-\frac{\hat{J}_{44}}{\lambda_4}\right) L} \\
e^{\left( \frac{\hat{J}_{22}}{v_2^*}-\frac{\hat{J}_{44}}{\lambda_4}\right) L} & 0 & 0 \\
0 & e^{\left( \frac{\hat{J}_{33}}{\lambda_3}-\frac{\hat{J}_{44}}{\lambda_4}\right) L} & 0 \\
\end{array}\right].
\end{align}
In addition, the input given in~\eqref{sec3:eq:ControlDesignModel_bcL} is defined as
\begin{equation}
\bar{U}(t) = e^{-\frac{\hat{J}_{44}}{\lambda_4}L}\hat{U}(t).
\end{equation}
All in all, the control design model is given by~\eqref{sec3:eq:ControlDesignModel_w123PDE}, \eqref{sec3:eq:ControlDesignModel_w4PDE},~\eqref{sec3:eq:ControlDesignModel_bc0} and~\eqref{sec3:eq:ControlDesignModel_bcL}. 
In Figure~\ref{sec3:fig:Qualitative_ControlDesignModel}, the qualitative behavior of the control design model is illustrated. According to the sign of the characteristic speeds, the propagation direction for each state $w_i(x,t)$ is drawn in Figure~\ref{sec3:fig:Qualitative_ControlDesignModel}. It shows that the control input $\bar{U}(t)$ acts at the outlet of system, first propagating upstream and, after it is carried through the boundary condition at the inlet of the investigated track section, affecting downstream traffic.

The summary of the two transformations is 
\begin{equation}
\label{sec3:eq:SummarizedTransformation}
\left[\begin{array}{c}
w_1 \\
w_2 \\
w_3 \\
w_4 \\
\end{array}\right] = T^{-1}(x)
\left[\begin{array}{c}
\tilde{\rho}_1 \\
\tilde{v}_1 \\
\tilde{\rho}_2 \\
\tilde{v}_2 \\
\end{array}\right]
\Leftrightarrow 
\left[\begin{array}{c}
\tilde{\rho}_1 \\
\tilde{v}_1 \\
\tilde{\rho}_2 \\
\tilde{v}_2 \\
\end{array}\right] = T(x)
\left[\begin{array}{c}
w_1 \\
w_2 \\
w_3 \\
w_4 \\
\end{array}\right],
\end{equation}
where
\begin{align}
&T^{-1}(x) = \left[\begin{array}{cccc}
0 & e^{-\frac{\hat{J}_{22}}{v_2^*}x} & 0 & 0 \\
0 & 0 & e^{-\frac{\hat{J}_{33}}{\lambda_3}x} & 0 \\
e^{-\frac{\hat{J}_{11}}{v_1^*}x} & 0 & 0 & 0 \\
0 & 0 & 0 & e^{-\frac{\hat{J}_{44}}{\lambda_4}x}
\end{array}\right]
V^{-1} , \label{sec3:eq:SummarizedTransformation_Tinv}\\
&T(x) = V\left[\begin{array}{cccc}
0 & 0 & e^{\frac{\hat{J}_{11}}{v_1^*}x} & 0 \\
e^{\frac{\hat{J}_{22}}{v_2^*}x} & 0 & 0 & 0 \\
0 & e^{\frac{\hat{J}_{33}}{\lambda_3}x} & 0 & 0 \\
0 & 0 & 0 & e^{\frac{\hat{J}_{44}}{\lambda_4}x} 
\end{array}\right]. \label{sec3:eq:SummarizedTransformation_T}
\end{align}
The overall transformation law depends on the spatial coordinate. In addition, both input transformations combined are
\begin{equation}
\label{sec3:eq:SummarizedTransformation_Ubar}
\bar{U}(t) = e^{-\frac{\hat{J}_{44}}{\lambda_4}L}\frac{1}{\kappa_4}U(t)
\end{equation}
and the inversion is given by
\begin{equation}
\label{sec3:eq:SummarizedTransformation_U}
U(t) = \kappa_4e^{\frac{\hat{J}_{44}}{\lambda_4}L}\bar{U}(t).
\end{equation}
Since all transformations are invertible, the stability properties of the linearized model in density and velocity perturbations and the UORM control design model are the same.
\section{Full-state feedback control design}
In the following, a full-state feedback control design for the system of four coupled hyperbolic PDEs given by~\eqref{sec3:eq:ControlDesignModel_w123PDE} and~\eqref{sec3:eq:ControlDesignModel_w4PDE} with boundary conditions~\eqref{sec3:eq:ControlDesignModel_bc0} and~\eqref{sec3:eq:ControlDesignModel_bcL} is carried out in order to achieve finite time convergence to zero for initial conditions $w_i(x,0) \in \mathcal{L}^\infty[0,L]$. The full-state feedback controller is designed by applying the backstepping control design in~\cite{HuKrsticGeneralPDEs}. The general idea is to transform the coupled hyperbolic PDEs to a cascade target system. The control law is chosen such that the instabilities in the system are eliminated through the boundary conditions of the target system. The states of the target system are denoted as $(\alpha_1,\alpha_2,\alpha_3,\beta)^T$. The kernels of the backstepping transformation are denoted by $K(x,\xi)$ and $L_{11}(x,\xi)$. Then, the backstepping transformation is defined as
\begin{align}
\alpha_1&(x,t) = w_1(x,t), \label{sec4:eq:Cntrl_BacksteppingTrafo1} \\
\alpha_2&(x,t) = w_2(x,t), \label{sec4:eq:Cntrl_BacksteppingTrafo2} \\
\alpha_3&(x,t) = w_3(x,t), \label{sec4:eq:Cntrl_BacksteppingTrafo3} \\
\beta&(x,t) = w_4(x,t) \notag \\
&-\int_0^x\left(K(x,\xi) \left[\begin{array}{c}
w_1(\xi,t) \\
w_2(\xi,t) \\
w_3(\xi,t) \\
\end{array}\right] + L_{11}(x,\xi)w_4(\xi,t)\right)d\xi \label{sec4:eq:Cntrl_BacksteppingTrafo4},
\end{align}
where
\begin{equation}
K(x,\xi) = \left[\begin{array}{ccc}
k_{11}(x,\xi) & k_{12}(x,\xi) & k_{13}(x,\xi)
\end{array}\right]
\end{equation}
and $L_{11}(x,\xi)$ are defined on a triangular domain
\begin{equation}
\label{sec4:eq:TriangularDomain_def}
\mathcal{T} = \{0\leq \xi\leq x\leq 1\}.
\end{equation}
The introduced kernels $K(x,\xi)$ and $L_{11}(x,\xi)$ are unknown and will be determined later on. Furthermore, the choice of the target system is 
\begin{align}
\left[\begin{array}{c}
\alpha_{1t}\\
\alpha_{2t} \\
\alpha_{3t}\\
\end{array} \right] =& -
\Lambda^+\left[\begin{array}{c}
\alpha_{1x} \\
\alpha_{2x} \\
\alpha_{3x} \\
\end{array}\right] +\Sigma^{++}(x)
\left[\begin{array}{c}
\alpha_1 \\
\alpha_2 \\
\alpha_3 \\
\end{array}\right] \notag \\
&+ \Sigma^{+-}(x)
\beta+ \int_0^x C^+(x,\xi)\left[\begin{array}{c}
\alpha_1(\xi,t) \\
\alpha_2(\xi,t) \\
\alpha_3(\xi,t) \\
\end{array}\right]d\xi \notag \\
&+ \int_0^x C^-(x,\xi)
\beta(\xi,t) d\xi, \label{sec4:eq:Cntrl_TargetSystem_a123PDE} \\
\beta_t=& \Lambda^-\beta_x \label{sec4:eq:Cntrl_TargetSystem_bPDE}.
\end{align}
The coefficients $C^+(x,\xi) \in \mathbb{R}^{3\times 3}$ and $C^-(x,\xi)\in \mathbb{R}^{3 \times 1}$ are defined on the same triangular domain $\mathcal{T}$ and are determined later on. Besides, the boundary conditions of the target system are
\begin{align}
\left[\begin{array}{c}
\alpha_1(0,t) \\
\alpha_2(0,t) \\
\alpha_3(0,t) \\
\end{array}\right] &= \bar{Q}_0 \beta(0,t), \label{sec4:eq:Cntrl_TargetSystem_BC0} \\
\beta(L,t) &= 0. \label{sec4:eq:Cntrl_TargetSystem_BCL}
\end{align}
The target system~\eqref{sec4:eq:Cntrl_TargetSystem_a123PDE} to~\eqref{sec4:eq:Cntrl_TargetSystem_BCL} converges to its equilibrium  at zero
\begin{equation}
\alpha_{e,i}(x)\equiv \beta_e(x) \equiv 0, \quad i=1,2,3, \text{ } t\geq 0, \text{ } x\in [0,L]
\end{equation}
in the finite time
\begin{equation}
\label{sec4:eq:FiniteConvergenceTime_tF}
t_F = \frac{L}{v_2^*}+\frac{L}{-\lambda_4}.
\end{equation}
The proof is given in Lemma $3.1$ in~\cite{HuKrsticGeneralPDEs}. It remains to compute the kernels $K(x,\xi)$ and $L_{11}(x,\xi)$, coefficients $C^+(x,\xi)$ and $C^-(x,\xi)$ and the control input $\bar{U}(t)$ such that the transformation is completed and to show the existence of the kernels. Deriving~\eqref{sec4:eq:Cntrl_BacksteppingTrafo4} with respect to space and time, inserting the resulting derivatives and~\eqref{sec3:eq:ControlDesignModel_bc0} in~\eqref{sec4:eq:Cntrl_TargetSystem_bPDE} yields the kernel equations that determine $K(x,\xi)$ and $L_{11}(x,\xi)$ after partial integration. The kernel equations are given by four coupled first order hyperbolic PDEs as well as four boundary conditions
\begin{align}
-\Lambda^- K_x(x,\xi) +K_\xi(x,\xi)\Lambda^+ =&-K(x,\xi)\Sigma^{++}(\xi) \notag \\
&-L_{11}(x,\xi)\Sigma^{-+}(\xi),\label{sec4:eq:kerneleq_PDEK} \\
-\Lambda^-L_{11x}(x,\xi)-L_{11\xi}(x,\xi)\Lambda^- =&-K(x,\xi)\Sigma^{+-}(\xi), \label{sec4:eq:kerneleq_PDEL} \\
K(x,0)\Lambda^+\bar{Q}_0-L_{11}(x,0)\Lambda^- =&\text{ }0, \label{sec4:eq:kerneleq_BCx0} \\
K(x,x)\Lambda^++\Lambda^-K(x,x) =&-\Sigma^{-+}(x). \label{sec4:eq:kerneleq_BCxx}
\end{align}
Inserting expressions for $\Lambda^-$, $K(x,\xi)$, $\Lambda^+$, $\Sigma^{++}(\xi)$, $\Sigma^{-+}(\xi)$, $\Sigma^{+-}(\xi)$ and denoting the entries of $\bar{Q}_0=\left\{\bar{Q}_{0i1}\right\}_{1\leq i\leq 3}$, the kernel equations in matrix form become
\begin{equation}
\begin{aligned}
&\lambda_4\left[\begin{array}{c}
k_{11x} \\
k_{12x} \\
k_{13x} \\
L_{11x}
\end{array}\right]+\left[\begin{array}{cccc}
v_2^* & 0 & 0 & 0 \\
0 & \lambda_3 & 0 & 0 \\
0 & 0 & v_1^* & 0 \\
0 & 0 & 0 & \lambda_4 \\
\end{array}\right]
\left[\begin{array}{c}
k_{11\xi} \\
k_{12\xi} \\
k_{13\xi} \\
L_{11\xi} \\
\end{array}\right] = \\
&\left[\begin{array}{cccc}
0 & -\bar{J}_{21}(\xi) & -\bar{J}_{31}(\xi) & -\bar{J}_{41}(\xi) \\
-\bar{J}_{12}(\xi) & 0 & -\bar{J}_{32}(\xi) & -\bar{J}_{42}(\xi) \\
-\bar{J}_{13}(\xi) & -\bar{J}_{23}(\xi) & 0 & -\bar{J}_{43}(\xi) \\
-\bar{J}_{14}(\xi) & -\bar{J}_{24}(\xi) & -\bar{J}_{34}(\xi) & 0 \\
\end{array}\right]\left[\begin{array}{c}
k_{11} \\
k_{12} \\
k_{13} \\
L_{11} \\
\end{array}\right]
\end{aligned}
\end{equation}
with boundary condition at $\xi=0$,
\begin{equation}
\left[\begin{array}{cccc}
\bar{Q}_{011}v_2^* & \bar{Q}_{021}\lambda_3 & \bar{Q}_{031}v_1^* & \lambda_4
\end{array}\right]\left[\begin{array}{c}
k_{11}(x,0) \\
k_{12}(x,0) \\
k_{13}(x,0) \\
L_{11}(x,0) \\
\end{array}\right] = 0,
\end{equation}
and boundary conditions at $\xi = x$,
\begin{align}
k_{11}(x,x) &= \frac{\bar{J}_{41}(x)}{\lambda_4-v_2^*}, \\
k_{12}(x,x) &= \frac{\bar{J}_{42}(x)}{\lambda_4-\lambda_3}, \\
k_{13}(x,x) &= \frac{\bar{J}_{43}(x)}{\lambda_4-v_1^*}. 
\end{align}
As shown in Theorem $3.3$ of~\cite{HuKrsticGeneralPDEs}, the kernel equations~\eqref{sec4:eq:kerneleq_PDEK} to~\eqref{sec4:eq:kerneleq_BCxx} are a well-posed system of equations and thus there exist unique solutions $K(x,\xi)$ and $L_{11}(x,\xi)$ in $L^\infty(\mathcal{T})$. Moreover, solving the equations~\eqref{sec4:eq:kerneleq_PDEL} and~\eqref{sec4:eq:kerneleq_BCx0} with the method of characteristics yields
\begin{align}
L_{11}(x,\xi) =& -\frac{1}{\lambda_4}K(x-\xi,0)\Lambda^+\bar{Q}_0 \notag \\
+\int_0^{-\frac{\xi}{\lambda_4}}&K(\lambda_4\nu+x,\lambda_4\nu+\xi)\Sigma^{+-}(\lambda_4\nu+\xi)d\nu. \label{sec4:eq:kerneleq_PDEL_Linserted}
\end{align}
Inserting this result in the remaining PDEs~\eqref{sec4:eq:kerneleq_PDEK} reduces the kernel equations to three coupled first order hyperbolic PDEs with three boundary conditions
\begin{align}
0=&\lambda_4K_x(x,\xi)+\Lambda^+K_\xi(x,\xi)+K(x,\xi)\Sigma^{++}(\xi) \notag \\
&-\frac{1}{\lambda_4}K(x-\xi,0)\Lambda^+\bar{Q}_0\Sigma^{-+}(\xi) \notag \\
+\int_0^{-\frac{\xi}{\lambda_4}}&K(\lambda_4\nu+x,\lambda_4 \nu+\xi)\Sigma^{+-}(\lambda_4\nu+\xi)d\nu\Sigma^{-+}(\xi) \label{sec4:eq:kerneleq_PDEK_Linserted} \\
0=&K(x,x)\Lambda^++\Lambda^-K(x,x)+\Sigma^{-+}(x) \label{sec4:eq:kerneleq_BCxx_Linserted}
\end{align}
Furthermore, deriving~\eqref{sec4:eq:Cntrl_BacksteppingTrafo1} to~\eqref{sec4:eq:Cntrl_BacksteppingTrafo3} with respect to space and time and inserting the obtained derivatives,~\eqref{sec4:eq:Cntrl_BacksteppingTrafo4} and~\eqref{sec3:eq:ControlDesignModel_w123PDE} in~\eqref{sec4:eq:Cntrl_TargetSystem_a123PDE} yields
\begin{align}
C^-(x,\xi) = \Sigma^{+-}(x)L(x,\xi)+\int_\xi^xC^-(x,s)L(s,\xi)ds, \\
C^+(x,\xi) = \Sigma^{+-}(x)K(x,\xi)+\int_\xi^xC^-(x,s)K(s,\xi)ds.
\end{align}
Finally, inserting~\eqref{sec3:eq:ControlDesignModel_bcL} and~\eqref{sec4:eq:Cntrl_TargetSystem_BCL} in~\eqref{sec4:eq:Cntrl_BacksteppingTrafo4} evaluated at $x=L$ delivers
\begin{align}
\bar{U}(t) &= -\bar{R}_1\left[\begin{array}{c}
w_1(L,t) \\
w_2(L,t) \\
w_3(L,t) \\
\end{array}\right] \notag \\
+&\int_0^L\left(K(L,\xi)\left[\begin{array}{c}
w_1(\xi,t) \\
w_2(\xi,t) \\
w_3(\xi,t) \\
\end{array}\right]+L_{11}(L,\xi) w_4(\xi,t)\right) d\xi
\end{align}
and therefore determines the control input. Before the results of the controller design are summarized in a theorem, the control law is formulated in dependence of the original physical variables, i.e. the densities  and velocities of both classes. For that reason, the transformation matrix $T^{-1}(x)$,~\eqref{sec3:eq:SummarizedTransformation_Tinv}, is separated in two parts
\begin{equation}
\label{sec3:eq:Seperation_Tinv}
T^{-1}(x) = \left[\begin{array}{c}
T_u^{-1}(x) \\
T_l^{-1}(x) \\
\end{array}\right],
\end{equation}
where $T_u^{-1}(x)\in \mathbb{R}^{3 \times 4}$ and $T_l^{-1}(x) \in \mathbb{R}^{1 \times 4}$. Hence, the states of the UORM control design model can be formulated as
\begin{align}
\left[\begin{array}{c}
w_1(\xi,t) \\
w_2(\xi,t) \\
w_3(\xi,t) \\
\end{array}\right] &= T_u^{-1}(\xi)\left[\begin{array}{c}
\tilde{\rho}_1(\xi,t) \\
\tilde{v}_1(\xi,t) \\
\tilde{\rho}_2(\xi,t) \\
\tilde{v}_2(\xi,t)
\end{array}\right], \\
w_4(L,t) &= T_l^{-1}(\xi)\left[\begin{array}{c}
\tilde{\rho}_1(\xi,t) \\
\tilde{v}_1(\xi,t) \\
\tilde{\rho}_2(\xi,t) \\
\tilde{v}_2(\xi,t)
\end{array}\right]
\end{align}
and the control law after applying the inverse input transformation~\eqref{sec3:eq:SummarizedTransformation_U} becomes
\begin{align}
&U(t) = -\kappa_4e^{\frac{\hat{J}_{44}}{\lambda_4}L}\bar{R}_1T_u^{-1}(L)\left[\begin{array}{c}
\rho_1(L,t)-\rho_1^* \\
v_1(L,t)-v_1^* \\
\rho_2(L,t)-\rho_2^* \\
v_2(L,t)-v_2^* \\
\end{array}\right] \notag \\
&-\kappa_4e^{\frac{\hat{J}_{44}}{\lambda_4}L}\int_0^LK(L,\xi)T_u^{-1}(\xi)\left[\begin{array}{c}
\rho_1(\xi,t)-\rho_1^* \\
v_1(\xi,t)-v_1^* \\
\rho_2(\xi,t) -\rho_2^* \\
v_2(\xi,t)-v_2^* \\
\end{array}\right]d\xi \notag \\
&-\kappa_4e^{\frac{\hat{J}_{44}}{\lambda_4}L}\int_0^LL_{11}(L,\xi)T_l^{-1}(\xi)\left[\begin{array}{c}
\rho_1(\xi,t)-\rho_1^* \\
v_1(\xi,t)-v_1^* \\
\rho_2(\xi,t) -\rho_2^* \\
v_2(\xi,t)-v_2^* \\
\end{array}\right]d\xi \label{sec4:eq:ControlLaw_PhysicalVariables}
\end{align}
using the definitions of the perturbations~\eqref{sec2:eq:PerturbationsDefinition1} and~\eqref{sec2:eq:PerturbationsDefinition2}. The result is now summarized in the following theorem.
\begin{thm}
Traffic density and velocity perturbations $(\tilde{\rho}_1,\tilde{\rho}_2,\tilde{v}_1,\tilde{v}_2)^T$ governed by the linearized two-class AR traffic model~\eqref{sec2:eq:LinearizedModelEq}, with the boundary conditions~\eqref{sec2:eq:BC0_lin_Dens1} to~\eqref{sec2:eq:BC0_lin_Flow} and~\eqref{sec2:eq:BCL_lin_Flow_wContrl} as well as initial profiles 
\begin{equation}
\tilde{\rho}_1(x,0),\tilde{v}_1(x,0),\tilde{\rho}_2(x,0),\tilde{v}_2(x,0) \in \mathcal{L}^\infty\left([0,L]\right),
\end{equation}
converge to the equilibrium at zero
\begin{equation}
\tilde{\rho}_{e,1}(x) \equiv \tilde{v}_{e,1}(x) \equiv \tilde{\rho}_{e,2}(x) \equiv \tilde{v}_{e,2}(x) \equiv 0
\end{equation}
in finite time $t_F$ given by~\eqref{sec4:eq:FiniteConvergenceTime_tF}, if the control law~\eqref{sec4:eq:ControlLaw_PhysicalVariables} is applied. Thereby, the kernels $K(x,\xi)$ and $L_{11}(x,\xi)$ are obtained by solving the kernel equations~\eqref{sec4:eq:kerneleq_PDEK_Linserted} and~\eqref{sec4:eq:kerneleq_BCxx_Linserted} on the triangular domain~\eqref{sec4:eq:TriangularDomain_def} and using~\eqref{sec4:eq:kerneleq_PDEL_Linserted} afterwards.
\end{thm}
The full-state feedback law in~\eqref{sec4:eq:ControlLaw_PhysicalVariables} requires measurements of the densities and velocities of both classes at every spatial point. In practice, it is possible to measure the densities and velocities at every spatial point by installing traffic cameras, collecting GPS data or helicopter data recordings.  
\section{Anti-collocated boundary observer design}
The installation of traffic cameras or the gain of GPS data in order to supply the full-state feedback control~\eqref{sec4:eq:ControlLaw_PhysicalVariables} with the required in-domain measurements is expensive. Therefore, a boundary observer design for full-state observation is proposed. In this work, an anti-collocated boundary observer is designed, i.e. the densities and velocities of both classes are measured at the opposite of the boundary where the control input acts. Therefore, it is assumed that only the traffic density and velocity of both classes at the inlet of the track section
\begin{align}
y_1(t) = \rho_1(0,t), \\
y_2(t) = v_1(0,t), \\
y_3(t) = \rho_2(0,t), \\
y_4(t) = v_2(0,t)
\end{align}
are measured. In terms of the control design model coordinates, inserting the measurements in~\eqref{sec3:eq:SummarizedTransformation} yields that
\begin{equation}
\bar{y}(t) = w_4(0,t)
\end{equation}
is known. 

The observer states $(\hat{w}_1,\hat{w}_2,\hat{w}_3,\hat{w}_4)^T$ are estimates of the control design model states $(w_1,w_2,w_3,w_4)^T$. Thus, the observer equations become
\begin{align}
\left[\begin{array}{c}
\hat{w}_{1t} \\
\hat{w}_{2t} \\
\hat{w}_{3t} \\
\end{array}\right] =& -\Lambda^+\left[\begin{array}{c}
\hat{w}_{1x} \\
\hat{w}_{2x} \\
\hat{w}_{3x} \\
\end{array}\right] + \Sigma^{++}(x)\left[\begin{array}{c}
\hat{w}_1 \\
\hat{w}_2 \\
\hat{w}_3 \\
\end{array}\right]+\Sigma^{+-}(x)\hat{w}_4 \notag \\
&-P^+(x)(\hat{w}_4(0,t)-w_4(0,t)) \label{sec5:eq:Observer_eqPDE_w1w2w3} \\
\hat{w}_{4t}=& \Lambda^-\hat{w}_{4x} + \Sigma^{-+}(x)\left[\begin{array}{c}
\hat{w}_1 \\
\hat{w}_2 \\
\hat{w}_3
\end{array}\right] \notag \\
&-P_{11}^-(x)(\hat{w}_4(0,t)-w_4(0,t)) \label{sec5:eq:Observer_eqPDE_w4}
\end{align}
where the gains of the output injections $P^+(x)$ and $P_{11}^-(x)$ need to be designed. The boundary conditions of the observer are
\begin{align}
\left[\begin{array}{c}
\hat{w}_1(0,t) \\
\hat{w}_2(0,t) \\
\hat{w}_3(0,t) \\
\end{array}\right] &= \bar{Q}_0w_4(0,t), \label{sec5:eq:Observer_eqbc0} \\
\hat{w}_4(L,t) &= \bar{R}_1\left[\begin{array}{c}
\hat{w}_1(L,t) \\
\hat{w}_2(L,t) \\
\hat{w}_3(L,t) \\
\end{array}\right] + \bar{U}(t). \label{sec5:eq:Observer_eqbcL}
\end{align}
As a next step, the system describing the dynamic behavior of the error between the states and their estimations is formulated. The estimation errors are defined as
\begin{equation}
\tilde{w}_i(x,t) = \hat{w}_i(x,t)-w_i(x,t), \quad i=1,2,3,4. 
\end{equation}
Subtracting the model equations of the control design model~\eqref{sec3:eq:ControlDesignModel_w123PDE},~\eqref{sec3:eq:ControlDesignModel_w4PDE},~\eqref{sec3:eq:ControlDesignModel_bc0} and~\eqref{sec3:eq:ControlDesignModel_bcL} from the observer equations~\eqref{sec5:eq:Observer_eqPDE_w1w2w3},~\eqref{sec5:eq:Observer_eqPDE_w4},~\eqref{sec5:eq:Observer_eqbc0} and~\eqref{sec5:eq:Observer_eqbcL} yields the following error system
\begin{align}
\left[\begin{array}{c}
\tilde{w}_{1t} \\
\tilde{w}_{2t} \\
\tilde{w}_{3t} \\
\end{array}\right] =&-\Lambda^+\left[\begin{array}{c}
\tilde{w}_{1x} \\
\tilde{w}_{2x} \\
\tilde{w}_{3x} \\
\end{array}\right] + \Sigma^{++}(x)\left[\begin{array}{c}
\tilde{w}_1 \\
\tilde{w}_2 \\
\tilde{w}_3 \\
\end{array}\right] \notag \\
&+ \Sigma^{+-}(x)\tilde{w}_4
-P^+(x)\tilde{w}_4(0,t) , \label{sec5:eq:ErrorSystem_w123PDE} \\
\tilde{w}_{4t} =& \Lambda^-\tilde{w}_{4x} + \Sigma^{-+}(x)\left[\begin{array}{c}
\tilde{w}_1 \\
\tilde{w}_2 \\
\tilde{w}_3 \\
\end{array}\right] \notag \\
&-P_{11}^-(x)\tilde{w}_4(0,t) \label{sec5:eq:ErrorSystem_w4PDE}
\end{align}
with the boundary conditions
\begin{align}
\left[\begin{array}{c}
\tilde{w}_1(0,t) \\
\tilde{w}_2(0,t) \\
\tilde{w}_3(0,t) \\
\end{array}\right] &=0, \label{sec5:eq:ErrorSystem_bc0}\\
\tilde{w}_4(L,t) &= \bar{R}_1\left[\begin{array}{c}
\tilde{w}_1(L,t) \\
\tilde{w}_2(L,t) \\
\tilde{w}_3(L,t) \\
\end{array}\right] \label{sec5:eq:ErrorSystem_bcL}.
\end{align}
Using the backstepping method, the output injection gains can be designed such that the error system converges to the equilibrium at zero in a finite time. Similar to the control design, a target system and a backstepping transformation are defined in the observer design as well. The output injections gains $P^+(x)$ and $P_{11}^-(x)$ are chosen such that the target converges to its equilibrium at zero in finite time. The state vector of the target system is denoted as $(\tilde{\alpha}_1,\tilde{\alpha}_2,\tilde{\alpha}_3,\tilde{\beta})^T$ and the kernels introduced in the backstepping transformation are $M(x,\xi)$ and $N_{11}(x,\xi)$. Thus, the backstepping transformation is given by
\begin{align}
\left[\begin{array}{c}
\tilde{w}_1(x,t) \\
\tilde{w}_2(x,t) \\
\tilde{w}_3(x,t) \\
\end{array}\right]  &= 
\left[\begin{array}{c}
\tilde{\alpha}_1(x,t) \\
\tilde{\alpha}_2(x,t) \\
\tilde{\alpha}_3(x,t) \\
\end{array} \right] + \int_0^xM(x,\xi)\tilde{\beta}(\xi,t) d\xi, \label{sec5:eq:Observ_BacksteppingTrafo_1}\\
\tilde{w}_4(x,t) &= \tilde{\beta}(x,t) + \int_0^xN_{11}(x,\xi)\tilde{\beta}(\xi,t)d\xi  \label{sec5:eq:Observ_BacksteppingTrafo_2}
\end{align}
where
\begin{equation}
M(x,\xi) = \left[\begin{array}{ccc}
m_{11}(x,\xi) & m_{21}(x,\xi) & m_{31}(x,\xi)
\end{array}\right]^T.
\end{equation}
The kernels $M(x,\xi)$ and $N_{11}(x,\xi)$ are defined in the triangular domain~\eqref{sec4:eq:TriangularDomain_def}. In addition, the target system is defined as
\begin{align}
\left[\begin{array}{c}
\tilde{\alpha}_{1t} \\
\tilde{\alpha}_{2t} \\
\tilde{\alpha}_{3t} \\
\end{array}\right] =& -\Lambda^+\left[\begin{array}{c}
\tilde{\alpha}_{1x} \\
\tilde{\alpha}_{2x} \\
\tilde{\alpha}_{3x} \\
\end{array} \right]+ \Sigma^{++}(x)\left[\begin{array}{c}
\tilde{\alpha}_1 \\
\tilde{\alpha}_2 \\
\tilde{\alpha}_3 \\
\end{array}\right] \notag \\
&+\int_0^xD^+(x,\xi)\left[\begin{array}{c}
\tilde{\alpha}_1(\xi,t) \\
\tilde{\alpha}_2(\xi,t) \\
\tilde{\alpha}_3(\xi,t) \\
\end{array}\right] d\xi, \label{sec5:eq:Observ_TargetSystem_alpha123PDE}\\
\tilde{\beta}_t =& \Lambda^-\tilde{\beta}_x + \Sigma^{-+}(x)\left[\begin{array}{c}
\tilde{\alpha}_1 \\
\tilde{\alpha}_2 \\
\tilde{\alpha}_3
\end{array}\right] \notag \\
&+\int_0^xD^-(x,\xi)\left[\begin{array}{c}
\tilde{\alpha}_1(\xi,t) \\
\tilde{\alpha}_2(\xi,t) \\
\tilde{\alpha}_3(\xi,t) \\
\end{array}\right] d\xi \label{sec5:eq:Observ_TargetSystem_betaPDE}
\end{align}
with the boundary conditions
\begin{align}
\left[\begin{array}{c}
\tilde{\alpha}_1(0,t) \\
\tilde{\alpha}_2(0,t) \\
\tilde{\alpha}_3(0,t) \\
\end{array}\right] &= 0, \label{sec5:eq:Observ_TargetSystem_bc0} \\
\tilde{\beta}(L,t) &= \bar{R}_1\left[\begin{array}{c}
\tilde{\alpha}_1(L,t) \\
\tilde{\alpha}_2(L,t) \\
\tilde{\alpha}_3(L,t) \\
\end{array}\right]. \label{sec5:eq:Observ_TargetSystem_bcL}
\end{align}
It can be shown that the target system converges in finite time $t_F$, given by~\eqref{sec4:eq:FiniteConvergenceTime_tF}. Besides, the coefficients $D^+(x,\xi) \in \mathbb{R}^{3 \times 3}$ and $D^-(x,\xi) \in \mathbb{R}^{1 \times 3}$ still need to be determined in the following. 

The equations for the output injection gains $P^+(x)$ and $P^-(x)$, the kernels of the backstepping transformation $M(x,\xi)$ and $N_{11}(x,\xi)$ and the coefficients in the target system $D^+(x,\xi)$ and $D^-(x,\xi)$ need to be deduced in a next step. The kernel equations for $M(x,\xi)$ and $N_{11}(x,\xi)$ are
\begin{align}
M_\xi(x,\xi)\Lambda^--\Lambda^+M_x(x,\xi) =&-\Sigma^{++}(x)M(x,\xi) \notag \\
&-\Sigma^{+-}(x)N_{11}(x,\xi), \label{sec5:eq:Observ_Kerneleq_MPDE} \\
N_{11\xi}(x,\xi)\Lambda^-+\Lambda^-N_{11x}(x,\xi)=&-\Sigma^{-+}(x)M(x,\xi), \label{sec5:eq:Observ_Kerneleq_NPDE} \\
M(\xi,\xi)\Lambda^-+\Lambda^+M(\xi,\xi)=& \Sigma^{+-}(\xi), \label{sec5:eq:Observ_Kerneleq_Mxx} \\
N_{11}(L,\xi)-\bar{R}_1M(L,\xi)=&0, \label{sec5:eq:Observ_Kerneleq_NLx}
\end{align}
where~\eqref{sec5:eq:Observ_Kerneleq_MPDE},~\eqref{sec5:eq:Observ_Kerneleq_NPDE} and~\eqref{sec5:eq:Observ_Kerneleq_Mxx} are obtained by inserting the transformation~\eqref{sec5:eq:Observ_BacksteppingTrafo_1} and~\eqref{sec5:eq:Observ_BacksteppingTrafo_2} as well as derivatives with respect to time and space of~\eqref{sec5:eq:Observ_BacksteppingTrafo_1} and~\eqref{sec5:eq:Observ_BacksteppingTrafo_2} in the PDEs of the error system~\eqref{sec5:eq:ErrorSystem_w123PDE} and~\eqref{sec5:eq:ErrorSystem_w4PDE}, followed by partial integration and noticing that $\tilde{\beta}(0,t)=\tilde{w}(0,t)$. In addition,~\eqref{sec5:eq:Observ_Kerneleq_NLx} is deduced by evaluating~\eqref{sec5:eq:Observ_BacksteppingTrafo_2} at $x=L$, plugging in the boundary conditions at the outlet of error system and target system,~\eqref{sec5:eq:ErrorSystem_bcL} and~\eqref{sec5:eq:Observ_TargetSystem_bcL}, and inserting~\eqref{sec5:eq:Observ_BacksteppingTrafo_1} evaluated at $x=L$ afterwards. Plugging in the expressions for $\Lambda^+$, $M(x,\xi)$, $\Lambda^-$, $\Sigma^{++}(x)$, $\Sigma^{+-}(x)$ and $\Sigma^{-+}(x)$, yields the kernel equations in matrix form:
\begin{align}
&\lambda_4\left[\begin{array}{c}
m_{11\xi}(x,\xi) \\
m_{21\xi}(x,\xi) \\
m_{31\xi}(x,\xi) \\
N_{11\xi}(x,\xi) \\
\end{array}\right] +\left[\begin{array}{cccc}
v_2^* & 0 & 0 & 0 \\
0 & \lambda_3 & 0 & 0 \\
0 & 0 & v_1^* & 0 \\
0 & 0 & 0 & \lambda_4 \\
\end{array}\right]\left[\begin{array}{c}
m_{11x}(x,\xi) \\
m_{21x}(x,\xi) \\
m_{31x}(x,\xi) \\
N_{11x}(x,\xi) \\
\end{array}\right] \notag \\
&= \left[\begin{array}{cccc}
0 & \bar{J}_{12}(x) & \bar{J}_{13}(x) & \bar{J}_{14}(x) \\
\bar{J}_{21}(x) & 0 & \bar{J}_{23}(x) & \bar{J}_{24}(x) \\
\bar{J}_{31}(x) & \bar{J}_{32}(x) & 0 & \bar{J}_{34}(x) \\
\bar{J}_{41}(x) & \bar{J}_{42}(x) & \bar{J}_{43}(x) & 0 \\
\end{array}\right]\left[\begin{array}{c}
m_{11}(x,\xi) \\
m_{21}(x,\xi) \\
m_{31}(x,\xi) \\
N_{11}(x,\xi) \\
\end{array}\right]
\end{align}
with boundary conditions at $x=\xi$ and $x=L$:
\begin{align}
m_{11}(\xi,\xi) &= \frac{\bar{J}_{14}(\xi)}{v_2^*-\lambda_4}, \\
m_{21}(\xi,\xi) &= \frac{\bar{J}_{24}(\xi)}{\lambda_3-\lambda_4}, \\
m_{31}(\xi,\xi) &= \frac{\bar{J}_{34}(\xi)}{v_1^*-\lambda_4}, \\
N_{11}(L,\xi) &= \bar{R}_1\left[\begin{array}{c}
m_{11}(L,\xi) \\
m_{21}(L,\xi) \\
m_{31}(L,\xi) \\
\end{array}\right].
\end{align}
It can be shown that the well-posedness of the kernel equations~\eqref{sec5:eq:Observ_Kerneleq_MPDE} to~\eqref{sec5:eq:Observ_Kerneleq_NLx} is equivalent to the kernel equations~\eqref{sec4:eq:kerneleq_PDEK} to~\eqref{sec4:eq:kerneleq_BCxx}. In fact, a transformation is introduced which achieves the exact same structure as the kernel equations which are developed during the full-state feedback design. Similar to the full-state feedback design, solving the PDE~\eqref{sec5:eq:Observ_Kerneleq_NPDE} and boundary condition~\eqref{sec5:eq:Observ_Kerneleq_NLx} with the method of characteristics delivers the expression
\begin{align}
N_{11}(x,\xi) =& \bar{R}_1M(L,L-(x-\xi))\notag \\
+\int_0^{\frac{x-L}{\lambda_4}}&\Sigma^{-+}(-\lambda_4\nu+x)M(-\lambda_4\nu+x,-\lambda_4\nu+\xi)d\nu, \label{sec5:eq:Observ_Kerneleq_NPDE_Ninserted}
\end{align}
in dependence of $M(x,\xi)$. Inserting this result in~\eqref{sec5:eq:Observ_Kerneleq_MPDE} reduces the kernel equations to three PDEs and three boundary conditions
\begin{align}
0=&-\Lambda^-M_\xi(x,\xi)+\Lambda^+M_x(x,\xi)-\Sigma^{++}(x)M(x,\xi)\notag \\
&-\Sigma^{+-}(x)\bar{R}_1M(L,L-(x-\xi))\notag \\
-\Sigma^{+-}&(x)\int_0^{\frac{x-L}{\lambda_4}}\Sigma^{-+}(x-\lambda_4\nu)M(x-\lambda_4\nu,\xi-\lambda_4\nu)d\nu, \label{sec5:eq:Observ_Kerneleq_MPDE_Ninserted} \\
0=&M(\xi,\xi)\Lambda^-+\Lambda^+M(\xi,\xi)-\Sigma^{+-}(\xi) \label{sec5:eq:Observ_Kerneleq_Mxx_Ninserted}
\end{align}
Besides, the computation that yields the kernel equations~\eqref{sec5:eq:Observ_Kerneleq_MPDE} to~\eqref{sec5:eq:Observ_Kerneleq_Mxx} for $M(x,\xi)$ and $N_{11}(x,\xi)$ implies
\begin{align}
D^+(x,\xi) =& -M(x,\xi)\Sigma^{-+}(\xi)+\int_\xi^xM(x,s)D^-(s,\xi)ds, \label{sec5:eq:Observ_D+} \\
D^-(x,\xi) =& -N_{11}(x,\xi)\Sigma^{-+}(\xi)\notag \\
&+\int_\xi^xN_{11}(x,s)D^-(s,\xi)ds, \label{sec5:eq:Observ_D-}
\end{align}
and
\begin{align}
P^+(x) &=-\lambda_4M(x,0), \label{sec5:eq:Observ_P+}\\
P_{11}^-(x) &= -\lambda_4N_{11}(x,0). \label{sec5:eq:Observ_P11-} 
\end{align}
Since the kernels $M(x,\xi)$ and $N_{11}(x,\xi)$ are well-posed,~\eqref{sec5:eq:Observ_D+} to~\eqref{sec5:eq:Observ_P11-} imply that the output injection gains as well as the target system coefficients are well-posed, too. Thus, the observer design is completed and is summarized in a theorem.
\begin{thm}
The error states $(\tilde{w}_1,\tilde{w}_2,\tilde{w}_3,\tilde{w}_4)^T$ between the observer~\eqref{sec5:eq:Observer_eqPDE_w1w2w3} to~\eqref{sec5:eq:Observer_eqbcL} and control design model~\eqref{sec3:eq:ControlDesignModel_w123PDE} to~\eqref{sec3:eq:ControlDesignModel_bcL} are described by~\eqref{sec5:eq:ErrorSystem_w123PDE} to~\eqref{sec5:eq:ErrorSystem_bcL}. If the output injections gains $P^+(x)$ and $P_{11}^-(x)$ are chosen as~\eqref{sec5:eq:Observ_P+} and~\eqref{sec5:eq:Observ_P11-}, where the kernel $M(x,\xi)$ is obtained by the well-posed equations~\eqref{sec5:eq:Observ_Kerneleq_MPDE_Ninserted} and~\eqref{sec5:eq:Observ_Kerneleq_Mxx_Ninserted} and $N_{11}(x,\xi)$ by~\eqref{sec5:eq:Observ_Kerneleq_NPDE_Ninserted}, and the initial error profiles are assumed to be
\begin{equation}
\tilde{w}_i (x,0) \in \mathcal{L}^\infty\left([0,L]\right), \quad i=1,2,3,4,
\end{equation}
then the errors converge to the equilibrium at zero
\begin{equation}
\tilde{w}_{e,i}(x)\equiv 0, \quad i=1,2,3,4
\end{equation}
in the finite time $t_F$ given by~\eqref{sec4:eq:FiniteConvergenceTime_tF}.
\end{thm}
The estimates of the observer $(\hat{w}_1,\hat{w}_2,\hat{w}_3,\hat{w}_4)^T$ can be transformed to the estimates of the density and velocity perturbations $(\hat{\tilde{\rho}}_1,\hat{\tilde{v}}_1,\hat{\tilde{\rho}}_2,\hat{\tilde{v}}_2)^T$ of both vehicle classes according to
\begin{equation}
\label{sec5:eq:SummarizedTransformation_estimtates}
\left[\begin{array}{c}
\hat{\tilde{\rho}}_1 \\
\hat{\tilde{v}}_1 \\
\hat{\tilde{\rho}}_2 \\
\hat{\tilde{v}}_2 \\
\end{array}\right] = T(x)
\left[\begin{array}{c}
\hat{w}_1 \\
\hat{w}_2 \\
\hat{w}_3 \\
\hat{w}_4 \\
\end{array}\right].
\end{equation}
Furthermore, the estimates of the original state variables are obtained by
\begin{align}
\label{sec5:eq:Perturbations_Def_estimates_dens}
\hat{\rho}_i(x,t) = \hat{\tilde{\rho}}_i(x,t)+\rho_i^*, \\
\hat{v}_i(x,t) = \hat{\tilde{v}}_i(x,t)+v_i^* \label{sec5:eq:Perturbations_Def_estimates_velo}
\end{align}
with respect to the estimates of the densities and velocities $(\hat{\rho}_1,\hat{v}_1,\hat{\rho}_2,\hat{v}_2)$.
\section{Output feedback control design}
So far, a full-state feedback, that requires measurements of all states at every spatial point and damps out stop-and-go traffic, and an observer, that generates estimates of all states at every spatial point based on a measurement at the inlet of the track section, has been designed. In a final step, both results are combined resulting in an output feedback control that damps out stop-and-go traffic based on a single measurement at the inlet of the track section. Therefore, the control law~\eqref{sec4:eq:ControlLaw_PhysicalVariables} is reformulated in terms of the generated estimates. This is done by replacing the densities and velocities by their estimates yielding the output feedback controller
\begin{align}
&U(t) = -\kappa_4e^{\frac{\hat{J}_{44}}{\lambda_4}L}\bar{R}_1T_u^{-1}(L)\left[\begin{array}{c}
\hat{\rho}_1(L,t)-\rho_1^* \\
\hat{v}_1(L,t)-v_1^* \\
\hat{\rho}_2(L,t)-\rho_2^* \\
\hat{v}_2(L,t)-v_2^* \\
\end{array}\right] \notag \\
&-\kappa_4e^{\frac{\hat{J}_{44}}{\lambda_4}L}\int_0^LK(L,\xi)T_u^{-1}(\xi)\left[\begin{array}{c}
\hat{\rho}_1(\xi,t)-\rho_1^* \\
\hat{v}_1(\xi,t)-v_1^* \\
\hat{\rho}_2(\xi,t) -\rho_2^* \\
\hat{v}_2(\xi,t)-v_2^* \\
\end{array}\right]d\xi \notag \\
&-\kappa_4e^{\frac{\hat{J}_{44}}{\lambda_4}L}\int_0^LL_{11}(L,\xi)T_l^{-1}(\xi)\left[\begin{array}{c}
\hat{\rho}_1(\xi,t)-\rho_1^* \\
\hat{v}_1(\xi,t)-v_1^* \\
\hat{\rho}_2(\xi,t) -\rho_2^* \\
\hat{v}_2(\xi,t)-v_2^* \\
\end{array}\right]d\xi \label{sec6:eq:OuputFeedbackControlLaw_PhysicalVariables}
\end{align}
where the estimates $(\hat{\rho}_1,\hat{v}_1,\hat{\rho}_2,\hat{v}_2)$ are obtained by transforming the state vector of the anti-collocated observer~\eqref{sec5:eq:Observer_eqPDE_w1w2w3} to~\eqref{sec5:eq:Observer_eqbcL} according to~\eqref{sec5:eq:SummarizedTransformation_estimtates},~\eqref{sec5:eq:Perturbations_Def_estimates_dens} and~\eqref{sec5:eq:Perturbations_Def_estimates_velo}, the transformation matrices $T_u^{-1}(\cdot)$ and $T_l^{-1}(\cdot)$ are given by~\eqref{sec3:eq:Seperation_Tinv}, the kernel $K(x,\xi)$ is the solutions of~\eqref{sec4:eq:kerneleq_PDEK_Linserted} and~\eqref{sec4:eq:kerneleq_BCxx_Linserted} and $L_{11}(x,\xi)$ is given by~\eqref{sec4:eq:kerneleq_PDEL_Linserted}. Finally, the abbreviation $\kappa_4$ is introduced in~\eqref{sec3:eq:kappa_abbrev}. This combination of the results obtained by the first two theorems is summarized in a final third theorem.
\begin{thm}
The linearized two-class AR model is given by~\eqref{sec2:eq:LinearizedModelEq} with the assumptions~\eqref{sec2:eq:BC0_lin_Dens1} to~\eqref{sec2:eq:BC0_lin_Flow} and~\eqref{sec2:eq:BCL_lin_Flow_wContrl} as boundary conditions. If the control law~\eqref{sec6:eq:OuputFeedbackControlLaw_PhysicalVariables} 
is applied in~\eqref{sec2:eq:BCL_lin_Flow_wContrl}, where the estimates are generated by the anti-collocated observer~\eqref{sec5:eq:Observer_eqPDE_w1w2w3} to~\eqref{sec5:eq:Observer_eqbcL} with the transformed control law~\eqref{sec3:eq:SummarizedTransformation_Ubar} as input and transformation~\eqref{sec5:eq:SummarizedTransformation_estimtates},~\eqref{sec5:eq:Perturbations_Def_estimates_dens} and~\eqref{sec5:eq:Perturbations_Def_estimates_velo} afterwards, and the initial profiles satisfy
\begin{equation}
\tilde{\rho}_1(x,0),\tilde{v}_1(x,0),\tilde{\rho}_2(x,0),\tilde{v}_2(x,0) \in \mathcal{L}^\infty\left([0,L]\right),
\end{equation}
then the perturbations converge to the equilibrium at zero
\begin{equation}
\tilde{\rho}_{e,1}(x) \equiv \tilde{v}_{e,1}(x) \equiv \tilde{\rho}_{e,2}(x) \equiv \tilde{v}_{e,2}(x) \equiv 0
\end{equation}
in the finite time $2t_F$, where $t_F$ is given by~\eqref{sec4:eq:FiniteConvergenceTime_tF}. The kernels $K(x,\xi)$ and $L_{11}(x,\xi)$ are obtained by solving the well-posed kernel equations~\eqref{sec4:eq:kerneleq_PDEK_Linserted},~\eqref{sec4:eq:kerneleq_BCxx_Linserted} and using~\eqref{sec4:eq:kerneleq_PDEL_Linserted} and the observer gains are given by~\eqref{sec5:eq:Observ_P+} and~\eqref{sec5:eq:Observ_P11-}, where the kernels $M(x,\xi)$ represents the solution of the well-posed system of equations~\eqref{sec5:eq:Observ_Kerneleq_MPDE_Ninserted} and~\eqref{sec5:eq:Observ_Kerneleq_Mxx_Ninserted} and $N_{11}(x,\xi)$ is given by~\eqref{sec5:eq:Observ_Kerneleq_NPDE_Ninserted}.
\end{thm}
As a final comment, the finite convergence time is twice as large as the time proposed in the previous theorems, because it requires $t_F$ to obtain correct state estimates and afterwards the control needs another $t_F$ to achieve convergence of the state variable to equilibrium state.
\section{Numerical simulation}
\begin{figure*}[htbp]
\begin{center}
% trim: links, unten, rechts, oben
\includegraphics[width=16cm,height=6cm, trim = {0cm 1cm 0cm 2cm},clip]{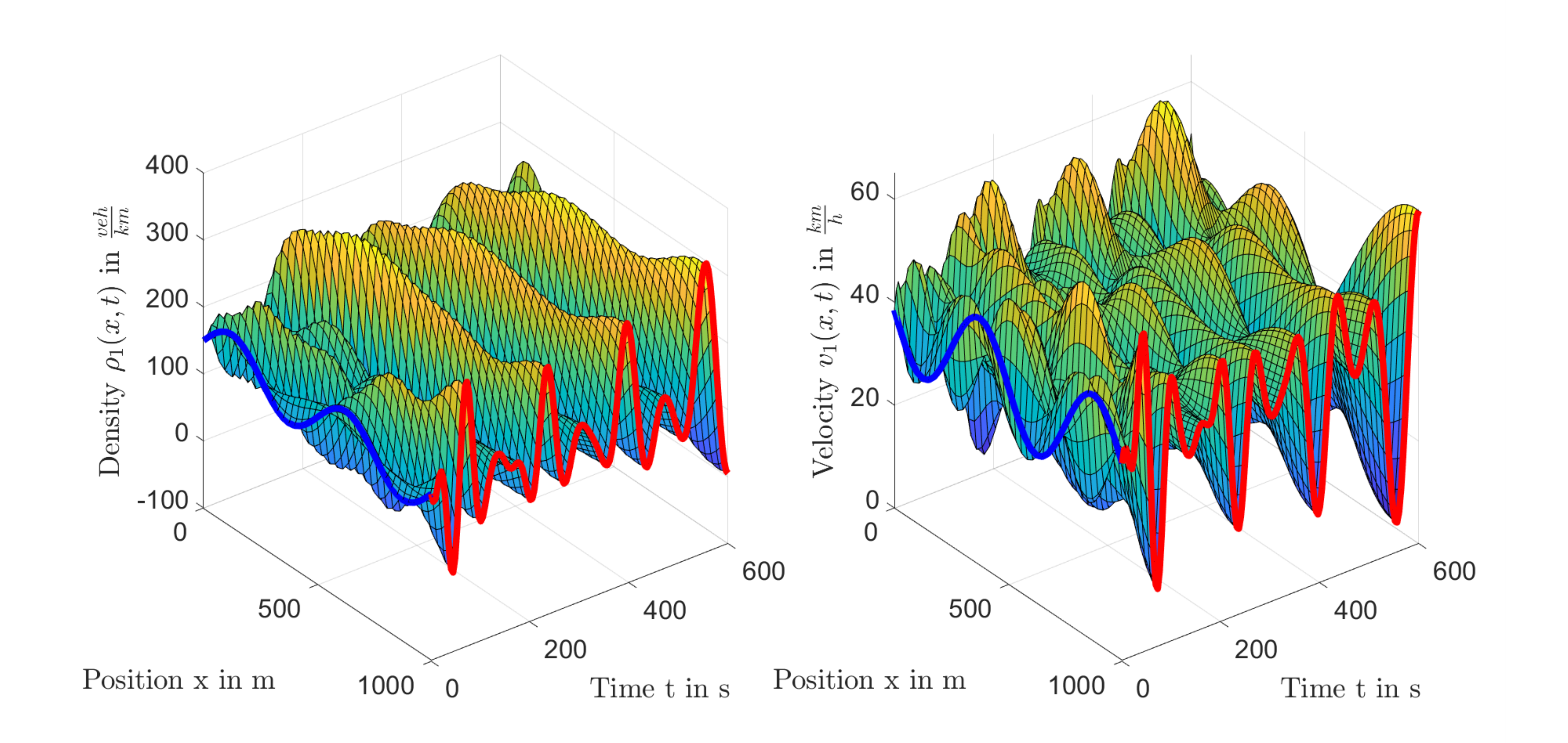}    % The printed column  
\caption{Traffic density and velocity of class $1$ without control.}	% width is 8.4 cm.
\label{sec6:fig:Class1_Openloop}                              % Size the figures 
\end{center}                                 % accordingly.
\end{figure*}

\begin{figure*}[htbp]
\begin{center}
% trim: links, unten, rechts, oben
\includegraphics[width=16cm,height=6cm, trim = {0cm 1cm 0cm 2cm},clip]{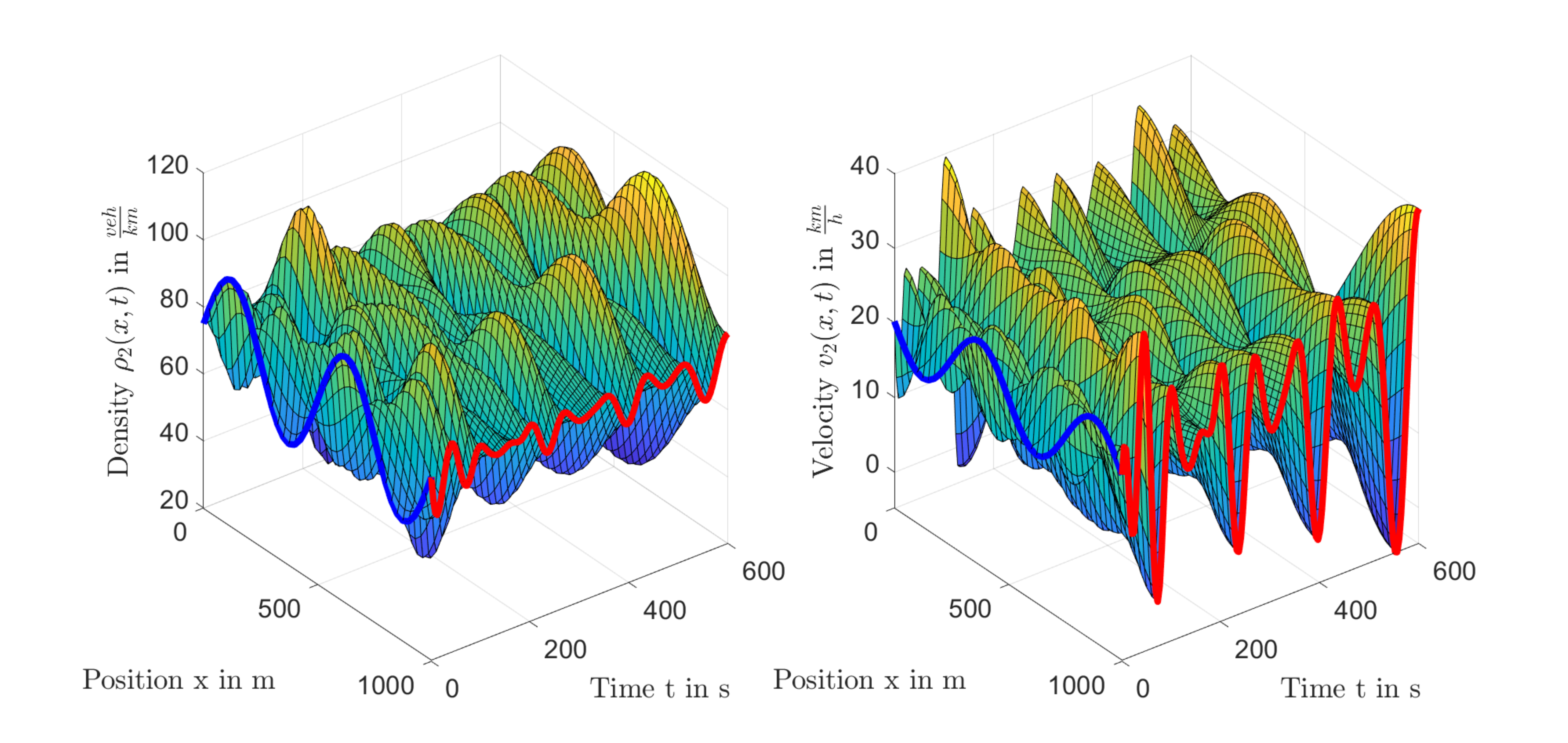}    % The printed column  
\caption{Traffic density and velocity of class $2$ without control.}	% width is 8.4 cm.
\label{sec6:fig:Class2_Openloop}                                % Size the figures 
\end{center}                                 % accordingly.
\end{figure*}

\begin{figure*}[htbp]
\begin{center}
% trim: links, unten, rechts, oben
\includegraphics[width=16cm,height=6cm, trim = {0cm 1cm 0cm 2cm},clip]{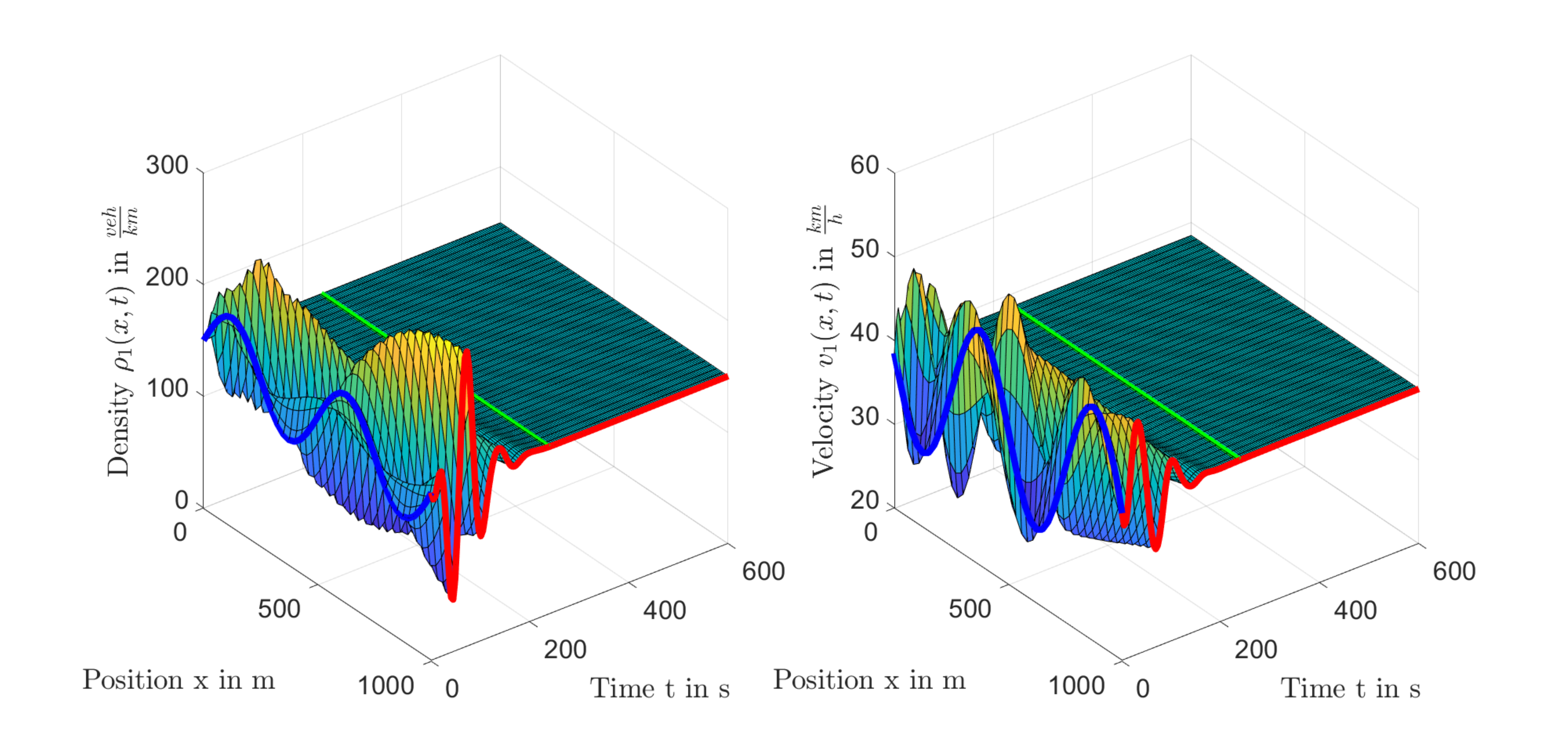}    % The printed column  
\caption{Traffic density and velocity of class $1$ with full-state feedback control. The green line indicates $t_F$.}	% width is 8.4 cm.
\label{sec6:fig:Class1_Control}                                 % Size the figures 
\end{center}                                 % accordingly.
\end{figure*}

\begin{figure*}[htbp]
\begin{center}
% trim: links, unten, rechts, oben
\includegraphics[width=16cm,height=6cm, trim = {0cm 1cm 0cm 2cm},clip]{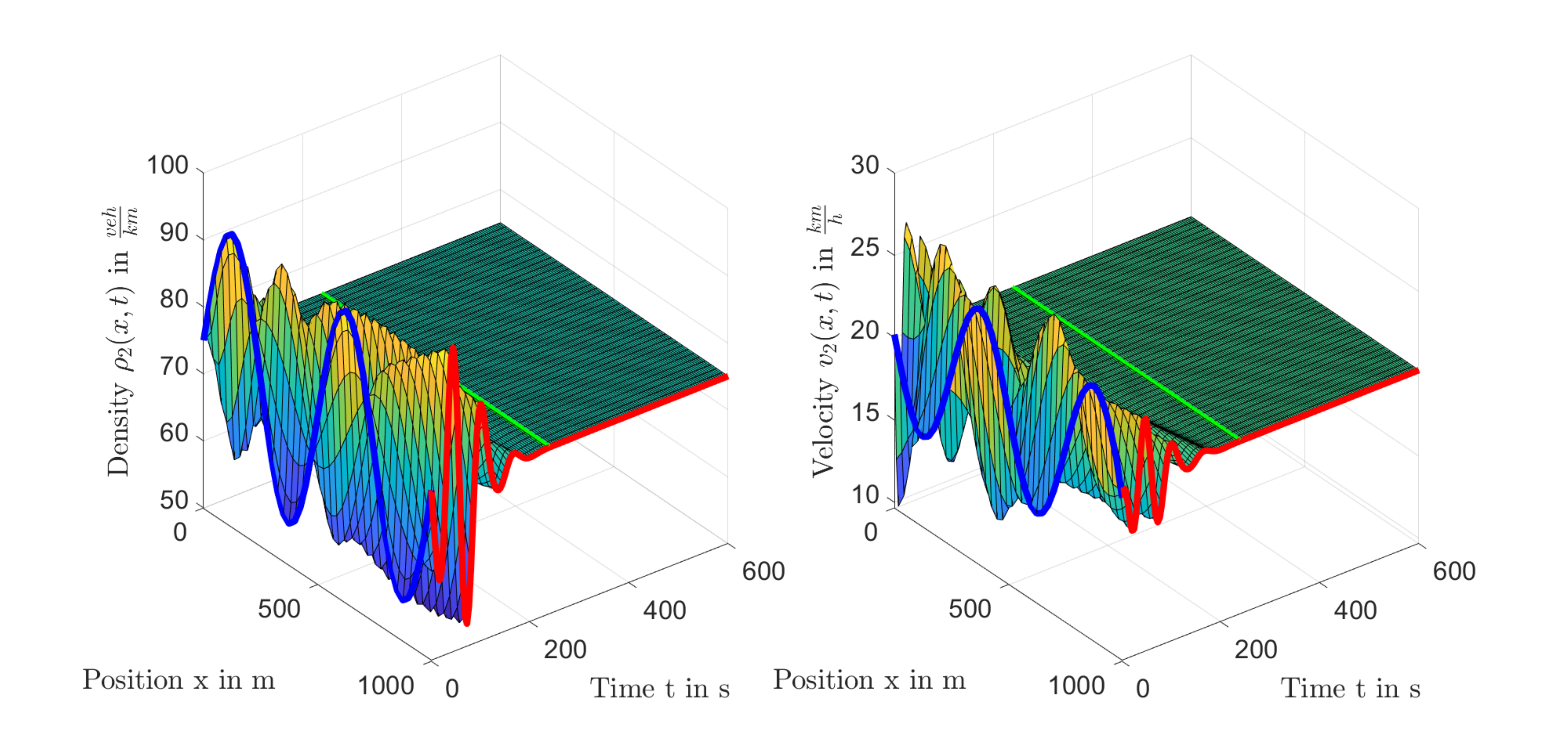}    % The printed column  
\caption{Traffic density and velocity of class $2$ with full-state feedback control. The green line indicates $t_F$.}	% width is 8.4 cm.
\label{sec6:fig:Class2_Control}                                 % Size the figures 
\end{center}                                 % accordingly.
\end{figure*}

\begin{figure*}[htbp]
\begin{center}
% trim: links, unten, rechts, oben
\includegraphics[width=16cm,height=6cm, trim = {0cm 1cm 0cm 2cm},clip]{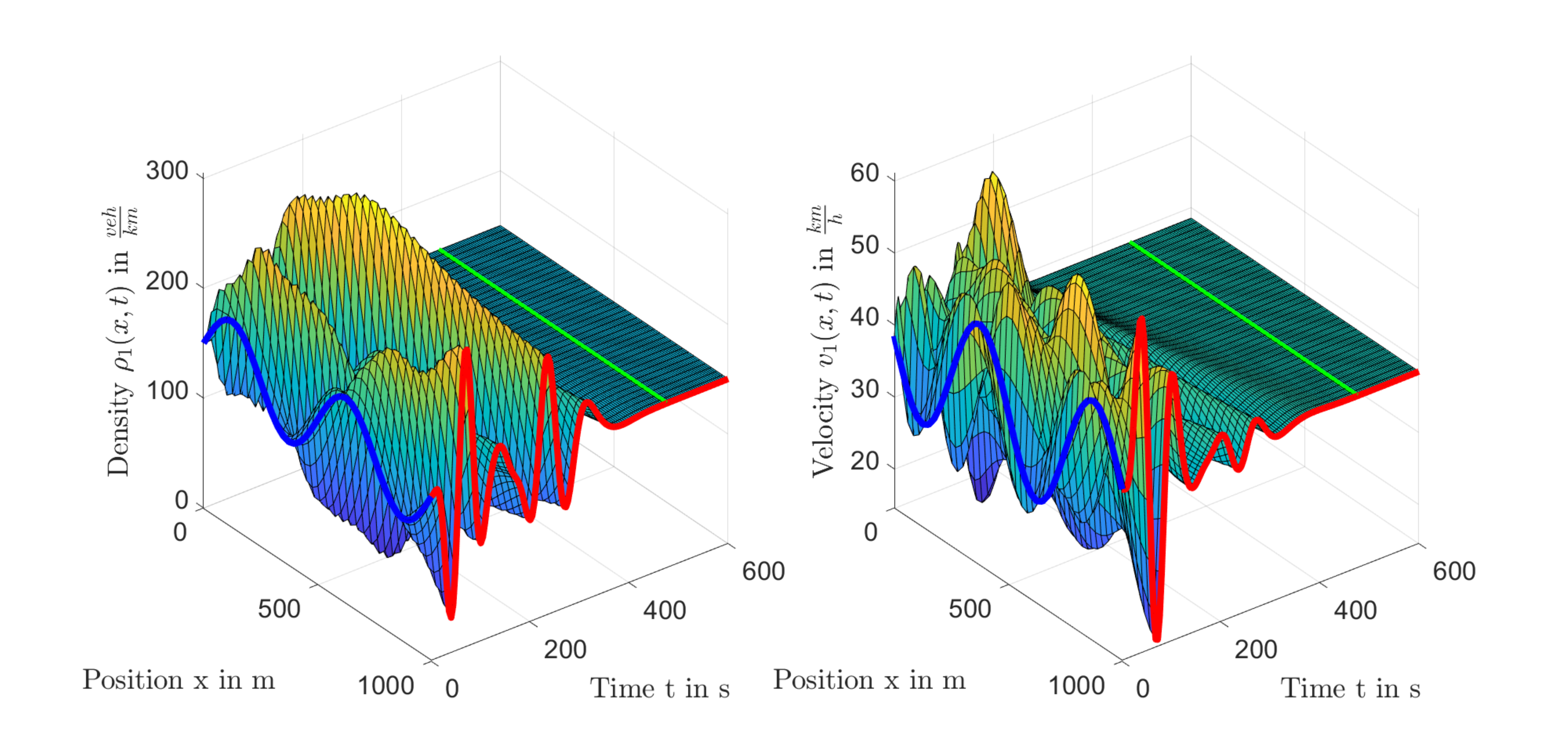}    % The printed column  
\caption{Traffic density and velocity of class $1$ with output feedback control. The green line indicates $2t_F$.}	% width is 8.4 cm.
\label{sec6:fig:Class1_Control_Observer}                                 % Size the figures 
\end{center}                                 % accordingly.
\end{figure*}

\begin{figure*}[htbp]
\begin{center}
% trim: links, unten, rechts, oben
\includegraphics[width=16cm,height=6cm, trim = {0cm 1cm 0cm 2cm},clip]{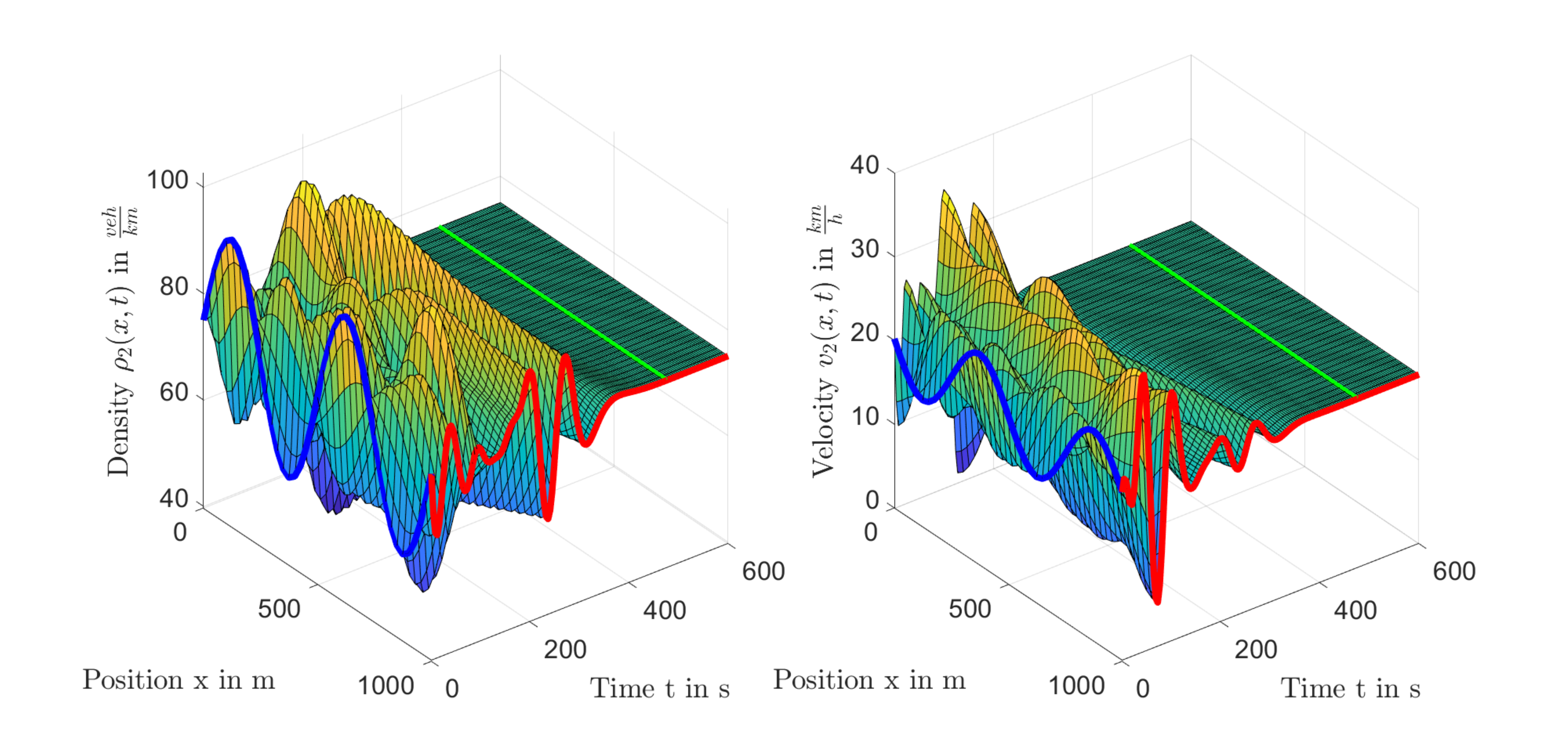}    % The printed column  
\caption{Traffic density and velocity of class $2$ with output feedback control. The green line indiciates $2t_F$.}	% width is 8.4 cm.
\label{sec6:fig:Class2_Control_Observer}                                 % Size the figures 
\end{center}                                 % accordingly.
\end{figure*}

In the end, the performance of the full-state feedback and output feedback control is investigated by simulation. The linearized model equations~\eqref{sec2:eq:LinearizedModelEq} are approximated by using an Upwind scheme. In order to achieve numerical stability, the grid sizes for the spatial coordinate and time are chosen such that the Courant-Friedrichs-Lewy condition,
\begin{equation}
\left|\frac{\lambda_i\Delta t}{\Delta x}\right| \leq 1, \quad i=1,2,3,4,
\end{equation}
is satisfied for all four characteristic speeds. The assumed parameter values are stated in Table~\ref{sec6:tab:Simulation_parameters}. In addition, the type of vehicles that are represented by the vehicle classes are denoted in Table~\ref{sec6:tab:Vehicle classes}. Typically, it holds that the larger the vehicle size the larger the relaxation time which is why $\tau_1<\tau_2$. Since vehicle class $1$ corresponds to smaller and faster average vehicles, the free-flow velocity $V_1$ is higher than $V_2$ and the highway needs to be occupied in a greater extent such that the average vehicles become jammed, i.e. $\overline{AO}_1>\bar{AO_2}$. Furthermore, it is assumed that faster and smaller vehicles experience less traffic pressure for low $AO$ values and therefore $\gamma_1>\gamma_2$. Finally, the equilibrium densities are chosen such that the investigated traffic is in the congested regime. The equilibrium velocities are determined by the choice of the equilibrium densities and result in $v_1^*\approx \unit[38]{\frac{km}{h}}$ and $v_2^* \approx \unit[20]{\frac{km}{h}}$. Although $v_1^*$ and $v_2^*$ seem to be low, the equilibrium velocities are realistic in case of congested traffic which is evenly distributed.
\begin{table}
\centering
\caption{Simulation parameters.}
\begin{tabular}{|c|c|c|c|}
\hline
Name & Symbol & Value & Unit \\
\hline
\hline
Relaxation time & $\tau_1$ & $30$ & $\unit{s}$ \\ 
 & $\tau_2$ & $60$ & $\unit{s}$ \\
 \hline
Pressure exponent & $\gamma_1$ & $2.5$ & $\unit{1}$ \\
 & $\gamma_2$ & $2$ & $\unit{1}$ \\
 \hline
Free-flow velocity & $V_{1}$ & $80$ & $\frac{\unit{km}}{\unit{h}}$ \\
& $V_{2}$ & $60$ & $\frac{\unit{km}}{\unit{h}}$ \\
\hline
Maximum $AO$ & $\overline{AO}_1$ & $0.9$ & $1$ \\
 & $\overline{AO}_2$ & $0.85$ & $1$ \\ 
\hline
Occupied surface per vehicle & $a_1$ & $10$ & $\unit{m^2}$ \\
 & $a_2$ & $40$ & $\unit{m^2}$ \\
\hline
Equilibrium density & $\rho_1^*$ & $150$ & $\frac{\unit{veh}}{\unit{km}}$ \\
 & $\rho_2^*$ & $75$ & $\frac{\unit{veh}}{\unit{km}}$ \\
\hline
Track width & $W$ & $6.5$ & $\unit{m}$ \\
\hline
Track length & $L$ & $1000$ & $\unit{m}$ \\
\hline
Amount of grid points & $N$ & $40$ & $\unit{1}$ \\
\hline 
\end{tabular}
\label{sec6:tab:Simulation_parameters}
\end{table}
\begin{table}
\centering
\caption{Traffic classes with length and width of each vehicle.}
\begin{tabular}{|c|c|c|c|}
\hline
Name & Class number & Length & Width \\
\hline
\hline
Average vehicle & $1$ & $\unit[5]{m}$ & $\unit[2]{m}$ \\ 
Big trucks & $2$ & $\unit[10]{m}$ & $\unit[4]{m}$ \\
\hline
\end{tabular}
\label{sec6:tab:Vehicle classes}
\end{table}
The initial profiles represent stop-and-go traffic with oscillations in density and velocity of sinusoidal shape. The mean value of the oscillations are the equilibrium values. At spatial points where the densities of both classes are increased, their velocities are decreased and thus the profiles
\begin{align}
\label{sec6:eq:Simulation_InitialProfiles}
\rho_i(x,0) = \rho_i^*+\frac{\rho_i^*}{4}\sin\left(\frac{4\pi}{L}x\right), \quad i=1,2, \\
v_i(x,0) = v_i^*-\frac{v_i^*}{4}\sin\left(\frac{4\pi}{L}x\right), \quad i =1,2
\end{align}
are assumed as initial profiles. \\
The simulation results of the open loop simulation are illustrated in Figure~\ref{sec6:fig:Class1_Openloop} for vehicle class $1$ and in Figure~\ref{sec6:fig:Class2_Openloop} for vehicle class $2$. In each figure, the left plot shows the density of the corresponding vehicle class, whereas the plot on the right hand side illustrates the velocity. The values of the states at the outlet of the track section are marked with a red line, whereas the blue line emphasizes the initial profiles~\eqref{sec6:eq:Simulation_InitialProfiles}. The four plots indicate that the stop-and-go oscillations do not vanish without the influence of control. Next, the Figures~\ref{sec6:fig:Class1_Control} and~\ref{sec6:fig:Class2_Control} illustrate the simulation results for the same initial condition but with activated full-state feedback control. The green line marks the finite convergence time $t_F\approx \unit[237]{s}$. Thus, it is easy to see that the convergence to the constant equilibrium profile in $t_F$ is achieved. Finally, Figure~\ref{sec6:fig:Class1_Control_Observer} and Figure~\ref{sec6:fig:Class2_Control_Observer} show the simulation results for the initial profiles using the designed output feedback control. Since the observer requires $t_F$ to estimate the states without error and afterwards the controller needs $t_F$ to achieve finite time convergence, the total finite convergence time is now $2t_F\approx \unit[474]{s}$ and therefore green line is adjusted accordingly.
\section{Concluding remarks}
This work leads to further problems that may be explored in the future. First, it is typically preferred that the measurement for the observer is at the same spot where the control input acts on the system. Therefore, the design of the collocated observer is a result of great interest. In addition, the extended AR traffic model presented in~\cite{MCAR} is formulated for $n$ classes and there are results for $n+m$ heterodirectional behaving linear PDEs in the literature which enables the extension to more than two classes and hence even more realistic considerations. Especially the definition of the congestion boundary in case of three or more classes would be an interesting result. Finally, a combination of the results presented in this work with the results regarding two lanes~\cite{2Lane_Yu} is of interest for further research.
\begin{ack}                              
Mark Burkhardt acknowledges financial support through the Baden-W\"urttemberg Stipendium of the Baden-W\"urttemberg Stiftung. In addition, he would like to thank Professor Oliver Sawodny for organizing the collaboration with Huan Yu and Miroslav Krstic. Furthermore, he acknowledges Kevin Schmidt for his help within fruitful discussions while this work was carried out. 
\end{ack}

%\bibliographystyle{plain}        % Include this if you use bibtex 
%\bibliography{references}           % and a bib file to produce the 
                                 % bibliography (preferred). The
                                 % correct style is generated by
                                 % Elsevier at the time of printing.

%\appendix
%\section{}    % Each appendix must have a short title.
								        % Sections and subsections are supported  
                                        % in the appendices.                                     
\end{document}